\documentclass[12pt]{amsart}
\usepackage{amssymb,amscd,latexsym,amsxtra}
\usepackage[all]{xy}
\usepackage[margin=4cm]{geometry}

\newtheorem{lem}[subsection]{Lemma}
\newtheorem{thm}[subsection]{Theorem}
\newtheorem{prop}[subsection]{Proposition}
\newtheorem{cor}[subsection]{Corollary}
\theoremstyle{definition}
\newtheorem{defn}[subsection]{Definition}
\newtheorem{rmk}[subsection]{Remark}

\newcommand{\C}{\mathbb{C}}
\newcommand{\N}{\mathbb{N}}
\newcommand{\Z}{\mathbb{Z}}
\newcommand{\Q}{\mathbb{Q}}
\newcommand{\U}{\mathfrak{U}}
\newcommand{\GL}{\mathsf{GL}}
\newcommand{\SL}{\mathsf{SL}}
\newcommand{\T}{\mathsf{T}}
\newcommand{\gl}{\mathfrak{gl}}
\renewcommand{\sl}{\mathfrak{sl}}
\newcommand{\E}{\mathsf{E}}
\renewcommand{\v}{\mathsf{v}}
\newcommand{\cf}{\operatorname{cf}}
\newcommand{\Sym}{\mathfrak{S}}
\newcommand{\End}{\operatorname{End}}
\newcommand{\Hom}{\operatorname{Hom}}
\newcommand{\sgn}{\mathrm{sign}}
\renewcommand{\c}{\mathbf{c}}
\newcommand{\cc}{\tilde{\c}}
\renewcommand{\d}{\mathbf{d}}
\renewcommand{\AA}{\tilde{A}}
\newcommand{\M}{\mathcal{M}}
\renewcommand{\le}{\leqslant}
\renewcommand{\ge}{\geqslant}

\newcommand{\Br}{\mathfrak{B}}

\newcommand{\I}{\mathbf{I}}
\newcommand{\ch}{\operatorname{ch}}
\newcommand{\tH}{H^\prime}
\newcommand{\tK}{K^\prime}
\newcommand{\Kbar}{\overline{K}}

\numberwithin{equation}{subsection}
\allowdisplaybreaks

\begin{document}
\title{The rational Schur algebra}

\author{Richard Dipper}
\address{Mathematisches Institut B,
Universit\"{a}t Stuttgart,
Pfaffenwaldring 57,
Stuttgart, 70569, Germany}
\email{Richard.Dipper@mathematik.uni-stuttgart.de}

\author{Stephen Doty}
\address{Mathematics and Statistics, 
Loyola University Chicago, 
Chicago, Illinois 60626 U.S.A.}
\email{doty@math.luc.edu} \thanks{Partially supported by DFG
project DI 531/5-2 and NSA grant DOD MDA904-03-1-00.}  

\thanks{The authors thank R.M.~Green for useful discussions. We are
also grateful to the referee for a suggestion that greatly simplified
\S\ref{7}.}


\begin{abstract} 
We extend the family of classical Schur algebras in type $A$, which
determine the polynomial representation theory of general linear
groups over an infinite field, to a larger family, the rational Schur
algebras, which determine the rational representation theory of
general linear groups over an infinite field. This makes it possible
to study the rational representation theory of such general linear
groups directly through finite dimensional algebras. We show that
rational Schur algebras are quasihereditary over any field, and thus
have finite global dimension.

We obtain explicit cellular bases of a rational Schur algebra by a
descent from a certain ordinary Schur algebra.  We also obtain a
description, by generators and relations, of the rational Schur
algebras in characteristic zero. 
\end{abstract}

\maketitle

\section*{Introduction}\noindent
A theme of contemporary representation theory is to approach
representations of infinite groups (e.g.\ algebraic groups) and
related algebraic structures through finite-dimensional algebras. For
instance, the polynomial representation theory of general linear
groups over an infinite field has been profitably studied from this
viewpoint, through the Schur algebras. In that particular case, the
efforts of various researchers culminated in \cite[Theorem
2.4]{Erdmann}, which showed that the problem of computing
decomposition numbers for such general linear groups is equivalent to
the same problem for symmetric groups. (Half of that equivalence was
known much earlier; see \cite[Theorem 3.4]{James}.)

The aim of this paper is to extend the classical Schur algebras
$S_K(n,r)$ to a larger family $S_K(n;r,s)$ of finite-dimensional
algebras, enjoying many similar properties as the Schur algebras,
which can be used to directly approach all the rational
representations (in the defining characteristic) of the general linear
groups over an infinite field, in the same way that the classical
Schur algebras approach the polynomial representations.
 
The paper is organized as follows. In \S\ref{1} of the paper, some
general lemmas, needed later, are proved. The reader should
probably skip \S1 and refer back as needed.  We formulate the main
definitions and formalisms, along with some basic results, in
\S\S\ref{2}--\ref{3}. These two sections may be regarded as an
extension of Green's book \cite{Green:book} to the more general
context considered here.  One new feature is the emergence of an
infinite dimensional algebra $S_K(n)_z$ which is useful for grading
the rational representations of the general linear group $\GL_n(K)$
($K$ an infinite field).  These infinite dimensional algebras are
obtained as an inverse limit of rational Schur algebras, and they may
be of independent interest.  Moreover, each rational Schur algebra
$S_K(n;r,s)$ is a quotient of $S_K(n)_{r-s}$.  For fixed $n\ge 2$, as
$r,s$ vary over the set $\N$ of non-negative integers, the family of
$S_K(n;r,s)$-modules is precisely the family of rational
$\GL_n(K)$-modules. 

In \S\ref{5} we give some combinatorial descriptions of the set $\pi =
\Lambda^+(n;r,s)$ of dominant weights defining each rational Schur
algebra.  In \S\ref{4} we show that the rational Schur algebras are
generalized Schur algebras, in the sense of \cite{Donkin:SA1}. It
follows that rational Schur algebras are always quasihereditary when
taken over a field, so have finite global dimension.  In \S\ref{5a} we
describe how to obtain cellular bases for rational Schur algebras, by
a descent from an ordinary Schur algebra.  In \S\ref{7} we formulate a
presentation of $S_K(n;r,s)$ by generators and relations (when $K$ has
characteristic zero) along the same lines as the earlier result
\cite{DG:PSA} for $S_K(n,r)$.

In \S\ref{6} we consider Schur--Weyl duality and describe the action of
a certain algebra $\Br_{r,s}^{(n)}$ centralizing (in characteristic
zero) the action of the general linear group $\GL_n(K)$ on mixed
tensor space ${\E_K}^{\otimes r} \otimes {\E_K^*}^{\otimes s}$. Thus
we can regard $S_K(n;r,s)$ as the centralizer algebra
\[
    S_K(n;r,s) = \End_{\Br_{r,s}^{(n)}}({\E_K}^{\otimes r} \otimes
    {\E_K^*}^{\otimes s})
\]
at least when $n \ge r+s$ and $K$ has characteristic zero.  The
restrictions (on characteristic and $n$) are believed unnecessary.

\section{General lemmas}\label{1}\noindent 
\subsection{}\label{gen:cf}
Let $\Gamma$ be any semigroup and $K$ an infinite field. Denote by
$K^\Gamma$ the $K$-algebra of $K$-valued functions on $\Gamma$, with
product $ff'$ of elements $f,f' \in K^\Gamma$ given by $s \to
f(s)\,f'(s)$ for $s\in \Gamma$. Given a representation $\tau: \Gamma
\to \End_K(V)$ in a $K$ vector space $V$, the {\em coefficient space}
of the representation is the subspace $\cf_\Gamma V$ of $K^\Gamma$
spanned by the coefficients $\{ r_{ab} \}$ of the representation.  The
coefficients $r_{ab} \in K^\Gamma$ are determined relative to a choice
of basis $v_a$ ($a \in I$) for $V$ by
\begin{equation}\label{gen:a}
\tau(g)\, v_b = \sum_{a \in I} r_{ab}(g)\, v_a
\end{equation}
for $g \in \Gamma$, $b \in I$. The {\em envelope} of the
representation $\tau$ is the subalgebra of $\End_K(V)$ generated by
the image of $\tau$. As we shall show below, the notions of
coefficient space and envelope, associated to a given representation
$\tau$, are dual to one another, at least when $V$ is of finite
dimension.

To formulate the result, let $K\Gamma$ be the semigroup algebra of
$\Gamma$.  Elements of $K\Gamma$ are sums of the form $\sum_{g\in
\Gamma} a_g g$ ($a_g \in K$) with finitely many $a_g \ne 0$. The group
multiplication extends by linearity to $K\Gamma$.  Note that $K\Gamma$
is also a coalgebra, with comultiplication given on generators by $x
\to x\otimes x$, and counit by $x \to 1$, for $x \in \Gamma$.  The
given representation $\tau: \Gamma \to \End_K(V)$ extends by linearity
to an algebra homomorphism $K\Gamma \to \End_K(V)$; denote this
extended map also by $\tau$. Obviously the envelope $[V]_\Gamma$ is
simply the image $\tau(K\Gamma)$. In other words, the representation
$\tau$ factors through its envelope: there is a commutative diagram
\begin{equation}
\begin{gathered}
\xymatrix{
K\Gamma \ar[rr]^\tau \ar@{->>}[dr] && \End_K(V)\\
& [V]_\Gamma \ar@{^{(}->}[ur]
}
\end{gathered}
\end{equation}
in which the leftmost and rightmost diagonal arrows are a surjection
and injection, respectively. Taking linear duals, the above
commutative diagram induces another one:
\begin{equation}
\begin{gathered}
\xymatrix{
(K\Gamma)^*   && \End_K(V)^* \ar[ll]_{\tau^*} \ar@{->>}[dl] \\
& {[V]_\Gamma}^* \ar@{_{(}->}[ul]
}
\end{gathered}
\end{equation}
There is a natural isomorphism of algebras $(K\Gamma)^* \simeq
K^\Gamma$, given by restricting a linear $K$-valued map on $K\Gamma$
to $\Gamma$; its inverse is given by the process of linearly extending
a $K$-valued map on $\Gamma$ to $K\Gamma$.  Note that the algebra
structure on $(K\Gamma)^*$ comes from dualizing the coalgebra
structure on $K\Gamma$.  Now we are ready for the aforementioned
result.

\begin{lem}\label{gen:lem}
The coefficient space $\cf_\Gamma(V)$ may be identified with the image
of $\tau^*$, so there is an isomorphism of vector spaces
${[V]_\Gamma}^* \simeq \cf_\Gamma V$. Moreover, if $V$ has finite
dimension, then $[V]_\Gamma \simeq (\cf_\Gamma V)^*$, again a vector
space isomorphism.  In that case $(\dim_K V < \infty)$ the first
isomorphism is an isomorphism of coalgebras and the second is an
isomorphism of algebras.
\end{lem}

\begin{proof}
We prove the first claim. Relative to the basis $v_a$ ($a \in I$) the
algebra $\End_K(V)$ has basis $e_{ab}$ ($a,b \in I$), where $e_{ab}$
is the linear endomorphism of $V$ taking $v_b$ to $v_a$ and taking all
other $v_c$, for $c \ne b$, to 0. In terms of this notation, equation
\eqref{gen:a} is equivalent with the equality
\begin{equation}
\tau(g) = \sum_{a,b \in I} r_{ab}(g) \,e_{ab}.
\end{equation}
Let $e'_{ab}$ be the basis of $\End_K(V)^*$ dual to the basis
$e_{ab}$, so that $e'_{ab}$ is the linear functional on $\End_K(V)$
taking the value 1 on $e_{ab}$ and taking the value 0 on all other
$e_{cd}$. Then one checks that $\tau^*$ carries $e'_{ab}$ onto
$r_{ab}$.  This proves that $\cf_\Gamma(V)$ may be identified with the
image of $\tau^*$. From the remarks preceding the statement of the
lemma, the rest of the claims now follow immediately. Note that the
coalgebra structure on ${[V]_\Gamma}^*$ is that induced by dualizing
the algebra structure on $[V]_\Gamma$, and similarly, the algebra
structure on $(\cf_\Gamma V)^*$ is induced by dualizing the coalgebra
structure on $\cf_\Gamma V$. (See the discussion following
\ref{def:rsa} ahead for more details on the latter point.)
\end{proof}

\begin{rmk}
The first statement of the lemma shows, in particular, that the
coefficient space of a representation does not depend on a choice of
basis for the representation.
\end{rmk}

We include the following lemma since we need it later and we are
unaware of a reference.

\begin{lem}\label{gen:lem2}
(a) Suppose $\Gamma$ is a group.  For any finite dimensional
$K\Gamma$-module $V$ we have an isomorphism
$(\cf_\Gamma(V^*))^{\mathrm{opp}} \simeq \cf_\Gamma(V)$ as coalgebras.

(b) Given a coalgebra $C$, let $C^\mathrm{opp}$ be its opposite
coalgebra. Then $(C^\mathrm{opp})^* \simeq (C^*)^{\mathrm{opp}}$ as
algebras.  That is, the linear dual of the opposite coalgebra may be
identified with the opposite algebra of the linear dual.
\end{lem}

\begin{proof}
(a) The isomorphism is given by $\widetilde{r}_{ab} \to r_{ab}$ where
the $r_{ab} \in K^\Gamma$ are the coefficients of $V$.  The $r_{ab}$
satisfy the equations
\begin{equation} \textstyle
g u_b = \sum_{a \in I} r_{ab}(g) u_a
\end{equation}
for all $g \in \Gamma$, where $(u_b)$ ($b \in I$) is a chosen basis of
$V$. The $\widetilde{r}_{ab}$ are the coefficients of $V^*$ with
respect to the dual basis, satisfying $\widetilde{r}_{ab}(g) =
r_{ab}(g^{-1})$ for all $g \in \Gamma$.

The $r_{ab}$ span $\cf_\Gamma (V)$ and satisfy
\begin{equation} \textstyle
\Delta(r_{ab}) =  \sum_{c\in I} r_{ac} \otimes r_{cb}
\end{equation}
for each $a,b \in I$. A calculation shows that the
$\widetilde{r}_{ab}$ (which span $\cf_\Gamma(V^*)$) satisfy
\begin{equation} \textstyle
\Delta(\widetilde{r}_{ab}) = \sum_{c\in I} 
\widetilde{r}_{cb}\otimes \widetilde{r}_{ac}.
\end{equation}
The opposite coalgebra structure is the one in which the tensors on
the right hand side of the preceding equality are reversed in order.
This proves part (a) of the lemma.

(b) This is easily checked and left to the reader.  
\end{proof}

\section{The categories $\M_K(n;r,s)$, $\M_K(n)_z$} \label{2}\noindent
\subsection{}
Let $n\ge 2$ be an integer and $K$ an arbitrary infinite field.
Henceforth we take $\Gamma = \GL_n(K)$, the group of nonsingular $n
\times n$ matrices with entries in $K$. The vector space $K^\Gamma$ of
$K$-valued functions on $\Gamma$ is naturally a (commutative,
associative) $K$-algebra with product $ff'$ of elements $f,f' \in
K^\Gamma$ given by $s \to f(s)\,f'(s)$ for $s\in \Gamma$.  The group
multiplication $\Gamma \times \Gamma \to \Gamma$ and unit element $1
\to \Gamma$ induce maps
\begin{equation}
\begin{CD} K^\Gamma @>\Delta>> K^{\Gamma \times \Gamma},\end{CD}
\qquad
\begin{CD} K^\Gamma @>\varepsilon>> K \end{CD}
\end{equation}
given by
\begin{equation}
\Delta(f) = [(s,t) \to f(st),\ \  s,t \in \Gamma], 
\qquad
\varepsilon(f) = f(1)
\end{equation}
and one easily checks that both $\Delta$, $\varepsilon$ are
$K$-algebra maps.  There is a map (tensor products are always taken
over $K$ unless we indicate otherwise)
\begin{equation}
K^\Gamma \otimes K^\Gamma \to K^{\Gamma \times \Gamma}
\end{equation}
determined by the condition $f\otimes f' \to [(s,t) \to f(s)\,f'(t),
\text{ all } s,t \in \Gamma]$. This map is injective and we use it to
identify $K^\Gamma \otimes K^\Gamma$ with its image in $K^{\Gamma
\times \Gamma}$. Thus we regard $K^\Gamma \otimes K^\Gamma$ as a
subspace of $K^{\Gamma \times \Gamma}$.

\subsection{}
We will need certain elements of $K^\Gamma$. For $1 \le i,j \le n$,
let $\c_{ij} \in K^\Gamma$ be defined by $\c_{ij}(A) =$ the $(i,j)$
entry of the matrix $A \in \Gamma$. Similarly, let $\cc_{ij} \in
K^\Gamma$ be defined by $\cc_{ij}(A) = \c_{ij}(A^{-1})$ for $A \in
\Gamma$.  Let $\d \in K^\Gamma$ be given by $\d(A) = \det A$ for $A
\in \Gamma$. Note that $\d^{-1} \in K^\Gamma$; that is, $\d$ is an
invertible element of $K^\Gamma$.  The assumption that $K$ is infinite
ensures that the $\{ \c_{ij} \}$ are algebraically independent amongst
themselves; similarly for the $\{ \cc_{ij} \}$.  By elementary
calculations with the definitions one verifies:
\begin{gather}
\Delta(\c_{ij}) = \sum_k \c_{ik} \otimes \c_{kj}; \qquad 
\varepsilon(\c_{ij}) = \delta_{ij}; \label{cat:a}\\
\Delta(\cc_{ij}) = \sum_k \cc_{kj} \otimes \cc_{ik}; \qquad 
\varepsilon(\cc_{ij}) = \delta_{ij}; \label{cat:b} \\
\Delta(\d) = \d \otimes \d; \qquad \varepsilon(\d) = 1. \label{cat:c}
\end{gather}
Here $\delta$ is the usual Kronecker delta: $\delta_{ij}$ is 1 if
$i=j$ and 0 otherwise.  By applying $\c_{ij}$ to the usual expression
in terms of matrix coordinates for $gg^{-1} = 1 = g^{-1}g$ ($g \in
\Gamma$) one verifies that
\begin{equation} \label{cat:f}
\sum_k \c_{ik} \cc_{kj} = \delta_{ij} = \sum_k \cc_{ik} \c_{kj}.
\end{equation}
The following formal identities also hold:
\begin{gather}
\d = \sum_{\pi \in \Sym_n} \sgn(\pi)\, \c_{1,\pi(1)} \cdots
\c_{n,\pi(n)}; \label{cat:d}\\
\d^{-1} = \sum_{\pi \in \Sym_n} \sgn(\pi)\,
\cc_{1,\pi(1)} \cdots \cc_{n,\pi(n)} \label{cat:e}
\end{gather}
where $\sgn(\pi)$ is the signature of a permutation $\pi$ in the
symmetric group $\Sym_n$ on $n$ letters.  

\subsection{}
By Cramer's rule, each $\cc_{ij}$ is expressible as a product of
$\d^{-1}$ with a polynomial expression in the variables $\c_{ij}$ ($1
\le i,j \le n$). Thus, the subalgebra $\AA_K(n)$ of $K^\Gamma$
generated by all $\c_{ij}$ together with $\d^{-1}$ coincides with the
subalgebra of $K^\Gamma$ generated by the all $\c_{ij}$ and all
$\cc_{ij}$. From formulas \eqref{cat:a} and \eqref{cat:b} it follows
that $\Delta \AA_K(n) \subset \AA_K(n)\otimes \AA_K(n)$; hence
$\AA_K(n)$ is a bialgebra with comultiplication $\Delta$ and counit
$\varepsilon$.  Actually, $\AA_K(n)$ is a Hopf algebra with antipode
induced from the inverse map $\Gamma \to \Gamma$; by \eqref{cat:f} the
antipode interchanges $\c_{ij}$ and $\cc_{ij}$.  As Hopf algebras, we
identify $\AA_K(n)$ with the affine coordinate algebra $K[\Gamma]$.
Let $A_K(n)$ be the subalgebra of $\AA_K(n)$ generated by all
$\c_{ij}$; this is the subspace spanned by all monomials in the
$\c_{ij}$ and it we identify it with the algebra of polynomial
functions on $\Gamma$.  Clearly, $\AA_K(n)$ is the localization of
$A_K(n)$ at $\d$.

\subsection{}
Given nonnegative integers $r,s$ let $\AA_K(n; r,s)$ be the subspace
of $\AA_K(n)$ spanned by all products of the form $\prod_{i,\,j}
(\c_{ij})^{a_{ij}} \prod_{i,\,j} ({\cc_{ij}})^{b_{ij}}$ such that
$\sum a_{ij} = r$ and $\sum b_{ij} = s$.  Then $\AA_K(n; r,s)$ is for
each $r,s$ a sub-coalgebra of $\AA_K(n)$. (Use the fact that $\Delta$
is an algebra homomorphism.)

Clearly $\AA_K(n) = \sum_{r,s} \AA_K(n;r,s)$. Because of \eqref{cat:f}
this sum is not direct; in fact for any $r,s$ we have for each $1\le i
\le n$ inclusions
\begin{equation}\label{cat:g}
\AA_K(n;r,s) \subset \AA_K(n;r+1,s+1)   
\end{equation}
since by \eqref{cat:f} any $f \in \AA_K(n;r,s)$ satisfies $f = f \sum_k
\c_{ik} \cc_{ki}$ and $f = f \sum_k \cc_{ik} \c_{ki}$, and the
right hand side of each equality is a sum of members of
$\AA_K(n;r+1,s+1)$.

\subsection{} \label{Inrs}
Let $\I(n,r)$ denote the set of all multi-indices $(i_1, \dots, i_r)$
such that each $i_a$ lies within the interval $[1,n]$, for $a = 1,
\dots, r$. Set $\I(n;r,s)=\I(n,r)\times \I(n,s)$. Given a pair $(I,J)$ of
elements of $\I(n;r,s)$, say $I = ((i_1,\dots,i_r),(i'_1,\dots,i'_s))$,
$J = ((j_1,\dots,j_r),(j'_1,\dots,j'_s))$, we write
\begin{equation}
  \c_{I,J} = \c_{i_1,j_1}\cdots \c_{i_r,j_r}\, 
             \cc_{i'_1,j'_1}\cdots \cc_{i'_s,j'_s}
\end{equation}
which is an element of $\AA_K(n;r,s)$. In fact, the set of all such
elements $\c_{I,J}$ spans $\AA_K(n;,r,s)$ as $(I,J)$ vary over
$\I(n;r,s)\times \I(n;r,s)$.  When working with this spanning set, one
must take the following equality rule into account:
\begin{equation}
  \c_{I,J} = \c_{L,M} \text{ if } (I,J) \sim (L,M)
\end{equation}
where we define $(I,J) \sim (L,M)$ if there exists some
$(\sigma,\tau)$ in $\Sym_r \times \Sym_s$ such that
$I(\sigma,\tau)=L$, $J(\sigma,\tau)=M$.
Here $\Sym_r \times \Sym_s$ acts on the right on $\I(n;r,s) =
\I(n,r)\times \I(n,s)$ by 
\[
((i_1,\dots,i_r),(i'_1,\dots,i'_s)) (\sigma,\tau) = 
((i_{\sigma(1)},\dots,i_{\sigma(r)}),(i'_{\tau(1)},\dots,i'_{\tau(s)}))
\]
The equality rule above states that the symbols $\c_{I,J}$ are
constant on $\Sym_r \times \Sym_s$-orbits for the above action of
$\Sym_r \times \Sym_s$ on $\I(n;r,s) \times \I(n;r,s)$.

Even after taking this equality rule into account, the spanning set
described above is not a basis of $\AA_K(n;r,s)$, since the set is not
linearly independent. One could obtain a basis of `bideterminants' for
$\AA_K(n;r,s)$, by dualizing the procedure in \S\ref{5a} ahead.

\subsection{Bidegree and degree}
We let $\M_K(n;r,s)$ denote the category of finite dimensional
$K\Gamma$-modules $V$ whose coefficient space $\cf_\Gamma(V)$ (see
\ref{gen:cf}) lies in $\AA_K(n;r,s)$.  (The morphisms between two
objects $V,V'$ in $\M_K(n;r,s)$ are homomorphisms of
$K\Gamma$-modules.)  Objects in $\M_K(n;r,s)$ afford rational
representations of the algebraic group $\Gamma = \GL_n(K)$.  In fact,
the rational representations of $\Gamma$ are precisely the
$K\Gamma$-modules $V$ for which $\cf_\Gamma(V)$ is contained in
$\AA_K(n)$. We call objects of $\M_K(n;r,s)$ rational representations
of $\Gamma$ of {\em bidegree} $(r,s)$. Note however that the bidegree
of a rational representation is not well-defined, since if $V$ is an
object of $\M_K(n;r,s)$ then $V$ is also an object of
$\M_K(n;r+1,s+1)$. But it is easy to see that every finite dimensional
rational $K\Gamma$-module has a unique {\em minimal} bidegree; this is
the bidegree $(r,s)$ such that $V$ belongs to $\M_K(n;r,s)$ but $V$
does not belong to $\M_K(n; r-1, s-1)$.

Given an integer $z\ge 0$ set $\AA_K(n)_z = \bigcup_{t\ge0}
\AA_K(n;z+t,t)$ and set $\AA_K(n)_{-z} = \bigcup_{t\ge0} \AA_K(n;t,z+t)$.
Then there is a direct sum decomposition
\begin{equation} \textstyle
  \AA_K(n) = \bigoplus_{z\in \Z} \AA_K(n)_z
\end{equation}
and each summand in this decomposition is a sub-coalgebra of
$\AA_K(n)$.  A finite dimensional rational $K\Gamma$-module is said to
be of rational {\em degree} $z \in \Z$ if its coefficient space lies
in $\AA_K(n)_z$.  Note that if $V$ is a homogeneous polynomial module
of degree $r$ then its rational degree is also $r$. Thus the notion of
degree just defined for rational modules extends the corresponding
notion as usually defined for polynomial representations.

If $V$ is infinite dimensional such that $V$ is a union of finite
dimensional rational submodules each of degree $z$, then we say that 
$V$ has degree $z$. 

Let $\M_K(n)_z$ be the category of rational $K\Gamma$-modules $V$ of
(rational) degree $z$, for any $z\in \Z$. (Again morphisms are just
$K\Gamma$-module homomorphisms.) Obviously $\M_K(n;r,s)$ is a
subcategory of $\M_K(n)_{r-s}$ for any $r,s$.

The equality $\AA_K(n) = \sum_{r,s} \AA_K(n;r,s)$ shows that as $r,s$
vary, the categories $\M_K(n;r,s)$ taken together comprise all of the
rational representations of $\Gamma$. We have the following more
precise result (compare with \cite[(2.2c)]{Green:book}).

\begin{thm}\renewcommand{\theenumi}{\alph{enumi}}\label{cat:thm}
Let $V$ be a rational $K\Gamma$-module.  

\noindent(a) $V$ has a direct sum decomposition of the form
$\textstyle V = \bigoplus_{z \in \Z} V_z$, where $V_z$ is a rational
$K\Gamma$-submodule of $V$ of degree $z$.

\noindent(b) For each $z \in \Z$, $V_z$ is expressible as a union of
an ascending chain of rational $K\Gamma$-submodules as follows:
\begin{enumerate}\renewcommand{\theenumi}{\roman{enumi}}
\item If $z\ge 0$ then $V_z$ is the union of a chain of the form
$V_{z,0} \subset V_{z+1,1} \subset \cdots \subset V_{z+t,t} \subset
\cdots$ where each $V_{z+t,t}$ is a rational $K\Gamma$-module of
bidegree $(z+t,t)$.
\item If $z\le 0$ then $V_z$ is the union of a chain of the form
$V_{0,-z} \subset V_{1,-z+1} \subset \cdots \subset V_{t,-z+t} \subset
\cdots$ where each $V_{t,-z+t}$ is a rational $K\Gamma$-module of
bidegree $(t,-z+t)$.
\end{enumerate}
\end{thm}

\begin{proof}
(a) The decomposition in (a) follows from a general theorem on
comodules; see \cite[(1.6c)]{Green:LFR}.

(b) Let $V_{r,s}$ in case $r-s = z$ be the unique maximal
$K\Gamma$-submodule of $V_z$ such that $\cf_\Gamma V_{r,s} \subset
\AA_K(n;r,s)$. 
\end{proof}

In other words, part (a) says that every rational representation of
$\GL_n(K)$ is a direct sum of homogeneous ones, where we call objects
of $\M_K(n)_z$ ($z\in \Z$) {\em homogeneous} of (rational) degree $z$.

\section{The rational Schur algebra $S_K(n;r,s)$}\label{3}\noindent
Theorem \ref{cat:thm} shows that each indecomposable rational
$K\Gamma$-module $V$ is homogeneous; {\em i.e.}, $V \in \M_K(n)_z$ for
some $z \in \Z$. Thus one may as well confine one's attention to
homogeneous modules.  Set $S_K(n)_z = (\AA_K(n)_z)^*$. Since
$\AA_K(n)_z$ is a $K$-coalgebra, $S_K(n)_z$ is a $K$-algebra in the
standard way.  The algebra $S_K(n)_z$ is in general infinite
dimensional. Since our main goal is to approach the rational
representations of $\Gamma$ through finite dimensional algebras, we
make the following definition.

\begin{defn}\label{def:rsa}
The {\em rational Schur algebra}, denoted by $S_K(n;r,s)$, is the
$K$-algebra $\AA_K(n;r,s)^* = \mathrm{Hom}_K(\AA_K(n;r,s), K)$, with
multiplication induced from the coproduct $\Delta$ on $\AA_K(n;r,s)$.
\end{defn}

The product of two elements $\xi, \xi'$ of $S_K(n;r,s)$ is computed by
first taking the tensor product $\xi \otimes \xi'$ of the two maps.
This gives a linear map from $\AA_K(n;r,s) \otimes \AA_K(n;r,s)$ to $K
\otimes K$. By identifying $K \otimes K$ with $K$ we regard the map
$\xi \otimes \xi'$ as taking values in the field $K$; it is customary
to write $\xi \overline{\otimes} \xi'$ for this slightly altered
map. Finally, by composing with $\Delta$ we obtain a linear map from $
\AA_K(n;r,s)$ to $K$; {\em i.e.}, an element of $S_K(n;r,s)$. So the product
$\xi\xi'$ is given by
\begin{equation}
 \xi \xi' = (\xi \overline{\otimes} \xi')\circ \Delta .
\end{equation}
Products in $S_K(n)_z$ are computed in a similar way. 

Since $\AA_K(n;r,s)$ is defined by a finite spanning set, it and its
dual $S_K(n;r,s)$ are finite dimensional, for any $r,s$.  Note that
the restriction of $\varepsilon$ to $\AA_K(n;r,s)$ is the identity
element of $S_K(n;r,s)$.

For each $r,s$ the inclusions $\AA_K(n;r,s) \subset \AA_K(n)_{r-s}$
and \eqref{cat:g} induce surjective algebra maps
\begin{equation}
  S_K(n)_{r-s} \to S_K(n;r,s), \quad  S_K(n;r+1,s+1) \to S_K(n;r,s).
\end{equation}
In fact, for $z\in \Z$ it is clear $\AA_K(n)_z$ is the union of one of
the two chains
\begin{equation}
\begin{gathered}
\AA_K(n;z,0) \subset \AA_K(n;z+1,1) \subset \cdots\ , \\
\AA_K(n;0,-z) \subset \AA_K(n;1,1-z) \subset \cdots
\end{gathered}
\end{equation}
according as $z\ge0$ or $z\le0$. It follows that $S_K(n)_z$ is the
inverse limit of the corresponding dual chain of rational Schur
algebras
\begin{equation}
  S_K(n)_z \simeq 
\begin{cases}
\underset{\longleftarrow}{\lim}_t  S_K(n;z+t,t) & \text{if } z\ge 0; \\
\underset{\longleftarrow}{\lim}_t  S_K(n;t,t-z) & \text{if } z\le 0. 
\end{cases}
\end{equation}
Since the rational Schur algebras determine the algebras $S_K(n)_z$ in
this simple way, it makes sense to focus our attention on them instead
of the more complicated infinite dimensional $S_K(n)_z$.

\subsection{}
J.A.~Green \cite{Green:book} formulated the notion of a (classical)
Schur algebra $S_K(n,r)$ for $\Gamma=\GL_n(K)$, extending some results
of Schur's dissertation from characteristic zero to arbitrary
characteristic. By definition, $S_K(n,r)$ is the linear dual of a
coalgebra $A_K(n,r)$; here $A_K(n,r)$ is simply the span of all
monomials in the $\c_{ij}$ of total degree $r$.  Note that the
classical Schur algebras in Green's sense are included among the
rational Schur algebras: the classical Schur algebra $S_K(n,r)$ is
just the rational Schur algebra $S_K(n;r,0)$. This follows immediately
from the equality $A_K(n,r)=\AA_K(n;r,0)$.

\subsection{} 
Following \cite[2.4]{Green:book}, for each $g \in \Gamma$ let $e_g \in
S_K(n;r,s)$ be determined by $e_g(c)=c(g)$ for any $c\in
\AA_K(n;r,s)$. We have $e_g e_{g'} = e_{gg'}$ for any $g,g' \in
\Gamma$; moreover, $e_1 = \varepsilon$. Thus by extending the map $g
\to e_g$ linearly one obtains a map $e: K\Gamma \to S_K(n;r,s)$, a
homomorphism of $K$-algebras. 

Any $f \in K^\Gamma$ has a unique extension to a linear map $f:
K\Gamma \to K$.  As discussed in \S\ref{1}, this gives an
identification $K^\Gamma \simeq (K\Gamma)^*$.  With this
identification, the image of an element $a = \sum a_g\, g$ of
$K\Gamma$ under $e$ is simply evaluation at $a$: $e(a)$ takes $c$ to
$c(a)$ for all $c\in \AA_K(n;r,s)$.

\begin{prop}
(a) The map $e: K\Gamma \to S_K(n;r,s)$ is surjective. 

\noindent(b) Let $f \in K^\Gamma$. Then $f \in \AA_K(n;r,s)$ if and
only if $f( \ker e ) = 0$.

\noindent(c) Let $V$ be a finite dimensional $K\Gamma$-module. Then
$V$ belongs to $\M_K(n;r,s)$ if and only if $(\ker e)V = 0$. 
\end{prop}

\begin{proof}
See \cite[2.4b,c]{Green:book}. The arguments given there are also
valid in the current context. 
\end{proof}

The proposition shows that the category $\M_K(n;r,s)$ is equivalent to
the category of finite dimensional $S_K(n;r,s)$-modules. An object $V$
in either category is transformed into an object of the other, using
the rule:
\begin{equation}
  a v = e(a) v, \qquad (a \in K\Gamma, v\in V). 
\end{equation}
Since both actions determine the same algebra of linear
transformations on $V$, the concepts of submodule, module
homomorphism, etc., coincide in the two categories.

\subsection{} \label{rsa:E}
Let $\E_K$ be the vector space $K^n$, regarded as column vectors. The
group $\Gamma = \GL_n(K)$ acts on $\E_K$, on the left, by matrix
multiplication; this action makes $\E_K$ into a rational
$K\Gamma$-module. The linear dual space $\E_K^* = \Hom_K(\E_K,K)$ is
also a rational $K\Gamma$-module, with $g\in \Gamma$ acting via $(g f)
= [v \to f(g^{-1}v), v\in \E_K]$.  Now if $\v_i$ ($1\le i \le n$) is
the standard basis of $\E_K$, {\em i.e.}\ $\v_i$ is the vector with a
1 in the $i$th position and 0 elsewhere, and if $\v'_i$ ($1\le i \le
n$) is the basis dual to the $\v_i$, then we have
\begin{equation}
  g\v_j = \sum_i \c_{ij}(g)\v_i, \quad
  g\v'_i = \sum_j \cc_{ij}(g)\v'_j
\end{equation}
which shows that $\cf_\Gamma(\E_K)$ is the subspace of $K^\Gamma$
spanned by the $\c_{ij}$ and $\cf_\Gamma(\E_K^*)$ is the subspace of
$K^\Gamma$ spanned by the $\cc_{ij}$.  Thus $\cf_\Gamma(\E_K) =
\AA_K(n;1,0)$ and $\cf_\Gamma(\E_K^*) = \AA_K(n;0,1)$.

It is well known that {\em coefficient spaces are multiplicative}, in
the following sense: if $V,W$ are two representations of $\Gamma$,
then $\cf_\Gamma(V\otimes W) = (\cf_\Gamma V)(\cf_\Gamma W)$, where
the product is taken in the algebra $K^\Gamma$. From this and the
preceding remarks it follows immediately that
\begin{equation}
  \AA_K(n;r,s) = \cf_\Gamma({\E_K}^{\otimes r} \otimes
  {\E_K^*}^{\otimes s}).
\end{equation}
Now by Lemma \ref{gen:lem} it follows immediately that 
\begin{equation}
  S_K(n;r,s) \simeq [{\E_K}^{\otimes r} \otimes {\E_K^*}^{\otimes
  s}]_\Gamma;
\end{equation}
{\em i.e.}, the rational Schur algebra in bidegree $r,s$ may be identified
with the envelope of the representation of $\Gamma$ on mixed tensor
space ${\E_K}^{r,s} := {\E_K}^{\otimes r} \otimes {\E_K^*}^{\otimes
s}$. So {\em we may identify $S_K(n;r,s)$ with the image of the
representation}
\begin{equation} \label{rho-rep}
  \rho_K: K\Gamma \to \End_K( {\E_K}^{r,s} ).
\end{equation}

\subsection{}
Next we obtain an alternative description of the rational Schur
algebras, as a quotient of the hyperalgebra $\U_K$. This will
ultimately lead to a proof that rational Schur algebras are
generalized Schur algebras.

Set $\U_\Q = \U_\Q(\gl_n)$, the universal enveloping algebra of the
Lie algebra $\gl_n$, over the rational field $\Q$.  The space $\E_\Q$
is naturally a $\U_\Q$-module, with action induced from the
$\Gamma$-action (the action of $\gl_n$ is given by left matrix
multiplication if we view elements of $\E_\Q$ as column vectors).  The
dual space $\E_\Q^*$ is also a $\U_\Q$-module in the usual way, by
regarding it as the dual module for the Lie algebra $\gl_n$, with
action 
\[
 (x \cdot f)(v) = -f(x v), \qquad \text{all } x\in \gl_n, f\in \E^*_\Q,
v\in \E_\Q;
\]
(this is the action induced from the $\Gamma$-action), and so
${\E_\Q}^{r,s} = {\E_\Q}^{\otimes r} \otimes {\E_\Q^*}^{\otimes s}$ is
naturally a $\U_\Q$-module.

It follows from the Poincar\'{e}--Birkhoff--Witt theorem that $\U_\Q$
is generated as an algebra by elements (which we shall denote by the
same symbols) corresponding to the matrices $X_{ij}$ ($1\le i, j \le
n$) in $\gl_n$.  Here, for any $1\le i,j\le n$, $X_{ij} =
(\delta_{ik}\delta_{jl})_{1\le k,l\le n}$ is the $n\times n$ matrix
with a unique 1 in the $i$th row and $j$th column, and 0 elsewhere. We
set $H_i:= X_{ii}$ for each $i=1, \dots, n$.  We shall need the
Kostant $\Z$-form $\U_\Z$, the subring (with 1) of $\U_\Q$ generated
by all
\begin{equation}
\tfrac{X_{ij}^p}{p!} \ \ (1\le i \ne j\le n); \qquad
\tbinom{H_{i}}{q} \ \ (1\le i \le n) 
\end{equation}
for any integers $p,q \ge 0$.

Set $\E_\Z = \Z^n = \U_\Z \v_1$ (note that $\v_1 \in \E_\Q$ is a
highest weight vector). The $\Z$-module $\E_\Z$ is an admissible
lattice in $\E_\Q$; similarly its dual space $\E_\Z^* = \Hom_\Z(\E_\Z,
\Z)$ is an admissible lattice in $\E_\Q^*$. We have the equalities
$\E_\Z = \sum \Z \v_i$ and $\E_\Z^* = \sum \Z \v'_i$.  It is clear
that there are isomorphisms of vector spaces
\begin{equation}
  \E_K \simeq \E_\Z \otimes_\Z K; \quad \E^*_K \simeq \E^*_\Z \otimes_\Z K
\end{equation}
for any field $K$. This induces an isomorphism of vector spaces
\begin{equation}
  {\E_K}^{r,s} \simeq {\E_\Z}^{r,s} \otimes_\Z K 
\end{equation}
where ${\E_\Z}^{r,s} = {\E_\Z}^{\otimes r} \otimes_\Z {\E_\Z^*}^{\otimes
s}$ for any $r,s$.  For any field $K$, set $\U_K = \U_\Z \otimes_\Z
K$. We have actions of $\U_K$ on $\E_K$, $\E_K^*$, ${\E_K}^{r,s}$
induced from the action of $\U_\Z$ on ${\E_\Z}$, $\E_\Z^*$,
${\E_\Z}^{r,s}$ by change of base ring.  Let
\begin{equation}\label{phi_K}
  \varphi_K: \U_K \to \End_K( {\E_K}^{r,s} ) 
\end{equation}
be the representation affording the $\U_K$-module structure on
${\E_K}^{r,s}$.

\begin{prop} \label{prop:Chev} 
For any infinite field $K$, the rational Schur algebra $S_K(n;r,s)$
may be identified with the image of the representation $\varphi_K:
\U_K \to \End_K( {\E_K}^{r,s} )$.
\end{prop}

\begin{proof}
This is essentially the Chevalley group construction.  Each
$\frac{X_{ij}^m}{m!} \in \U_\Z$ ($1\le i\ne j \le n$) induces a
corresponding element $\frac{X_{ij}^m}{m!} \otimes 1$ of $\U_K$.  For
$t \in K$ and each $i \ne j$ let $x_{ij}(t)$ be the element
$I+tX_{ij}$ of $\SL_n(K)$, where $I$ is the $n\times n$ identity
matrix. The action of $\Gamma=\GL_n(K)$ on ${\E_K}^{r,s}$ is given by
the representation $\rho_K$ of \eqref{rho-rep}.  The element
$x_{ij}(t)$ acts on ${\E_K}^{r,s}$ as the $K$-linear endomorphism
$\rho_K(x_{ij}(t))$, where
\begin{equation}\label{Eij} \textstyle
\rho_K(x_{ij}(t)) = \sum_{m\ge 0} t^m
\varphi_K(\frac{X_{ij}^m}{m!}\otimes 1).
\end{equation}
The sum is finite because $X_{ij}$ acts nilpotently on ${\E_\Z}^{r,s}$.
Hence for each fixed $i \ne j$ there exists a natural number $N$ such
that
\begin{equation}\label{Eij-long} \textstyle
\rho_K(x_{ij}(t)) = 1 + t \varphi_K(X_{ij}\otimes 1) + \cdots + t^N
\varphi_K(\frac{X_{ij}^N}{N!})\otimes 1)
\end{equation}
for all $t \in K$.  Choosing $N+1$ distinct values $t_0$, $t_1, \dots,
t_N$ for $t$ in $K$ (which is always possible since $K$ is infinite)
we can solve the resulting linear system, by inverting the Vandermonde
matrix of coefficients, to obtain equalities
\begin{equation}\label{inverted} \textstyle
\varphi_K(\frac{X_{ij}^m}{m!} \otimes 1) = \sum_{k} a_{mk}\,
\rho_K(x_{ij}(t_k))
\end{equation}
for certain $a_{mk} \in K$.

The elements $x_{ij}(t)$ ($1\le i\ne j \le n$, $t\in K$) generate
the group $\SL_n(K)$. Equations \eqref{Eij} show that the image
$\rho_K(K\SL_n(K))$ is contained in $\varphi_K( \U'_K )$ where $\U'_K$
is the $K$-subalgebra of $\U_K$ generated by all $\frac{X_{ij}^m}{m!}
\otimes 1$. On the other hand, equations \eqref{inverted} justify the
opposite inclusion. Thus we have proved that 
\begin{equation} \label{SL-image}
  \rho_K(K\SL_n(K)) = \varphi_K( \U'_K ).
\end{equation}

For the algebraic closure $\Kbar$ of $K$ the group $\GL_n(\Kbar)$ is
generated by $\SL_n(\Kbar)$ together with the scalar matrices.  Since
the scalar matrices act as scalars on ${\E_{\Kbar}}^{r,s}$ we have the
equality $\rho_{\Kbar}( \Kbar\GL_n(\Kbar) ) =
\rho_{\Kbar}(\Kbar\SL_n(\Kbar))$, as the scalar operators are clearly
already present in $\rho_{\Kbar}(\Kbar\SL_n(\Kbar))$.  For general
(infinite) $K$ we have by construction equalities of dimension:
\begin{gather*}
\dim_K \rho_K(K\SL_n(K)) = \dim_{\Kbar}
\rho_{\Kbar}(\Kbar\SL_n(\Kbar)); \\
\dim_K \rho_K(K\GL_n(K)) = \dim_{\Kbar}
\rho_{\Kbar}(\Kbar\GL_n(\Kbar))
\end{gather*}
and hence the natural inclusion $\rho_K( K\SL_n(K) ) \subseteq
\rho_K(K\GL_n(K))$ must be an equality.

Moreover, $\U_K$ is generated by $\U'_K$ along with the elements
$\binom{H_i}{m} \otimes 1$.  By restricting $\varphi_K$ to the
isomorphic copy of $\U(\gl_2)$ generated by $X_{ij}$, $X_{ji}$, $H_i$,
and $H_j$ (for $i\ne j$) one can show by calculations similar to those
in \cite{DG:n=2} that the $\varphi_K(\binom{H_i}{m} \otimes 1)$ are
already present in $\varphi_K(\U'_K)$. (See also \ref{n=2case} ahead;
this justifies the applicability of \cite{DG:n=2}.) It follows that
$\varphi_K(\U'_K) = \varphi_K(\U_K)$, and the proof is complete.
\end{proof}

We note that the algebra $\U'_K$ appearing in the proof is simply the
hyperalgebra of $\sl_n$, obtained from the Kostant $\Z$-form of
$\U_\Q(\sl_n)$ by change of base ring. The above proof reveals also
the following.

\begin{cor} \label{SL-cor}
For any infinite field $K$, $S_K(n;r,s)$ is the image of the
restricted representation $K\SL_n(K) \to \End_K({\E_K}^{r,s})$.  This
may be identified with the image of the restricted representation
$\U'_K \to \End_K({\E_K}^{r,s})$.
\end{cor}

The inverse anti-automorphism $g \to g^{-1}$ on $\Gamma$ induces the
following result, which in particular says that $S_K(n;0,s)$ is
isomorphic with the opposite algebra of the ordinary Schur algebra
$S_K(n,s)$.

\begin{prop} \label{opp}
Let $K$ be an infinite field. For any $n,r,s$ we have an isomorphism
$S_K(n;r,s)^{\mathrm{opp}} \simeq S_K(n;s,r)$.
\end{prop}

\begin{proof}
Clearly we have an isomorphism $({\E_K}^{s,r})^* \simeq {\E_K}^{r,s}$,
as rational $K\Gamma$-modules. Applying Lemma \ref{gen:lem2}(a) to $V
= {\E_K}^{s,r}$ we obtain an isomorphism of coalgebras $(\cf_\Gamma
({\E_K}^{r,s}))^{\mathrm{opp}} \simeq \cf_\Gamma {\E_K}^{s,r}$.
Thus we have an isomorphism of coalgebras $\AA_K(n;r,s)^{\mathrm{opp}}
\simeq \AA_K(n;s,r)$, and the result follows by taking linear duals,
using Lemma \ref{gen:lem2}(b).
\end{proof}

\section{Weights of $\M_K(n;r,s)$}\label{5}\noindent
\subsection{} \label{wts}
Elements of the diagonal torus $\T_n = \{ \mathrm{diag}(t_1,\dots,t_n)\in
\GL_n(K)\}$ act semisimply on any object of the category
$\M_K(n;r,s)$, where the action is the one obtained from the action of
$\Gamma=\GL_n(K)$ by restriction.  Given an object $V$ in
$\M_K(n;r,s)$, in particular $V$ is a rational $K\Gamma$-module, so
$V$ is a direct sum of its weight spaces:
\begin{equation}
V = \bigoplus_{\lambda \in \Z^n} V_\lambda
\end{equation}
where $V_\lambda = \{ v\in V: \mathrm{diag}(t_1,\dots,t_n)\,v =
t_1^{\lambda_1}\cdots t_n^{\lambda_n} \,v \}$. Here we write
$\lambda$ for the vector $(\lambda_1, \dots, \lambda_n) \in \Z^n$.
The set of $\lambda$ for which $V_\lambda \ne 0$ is the set of weights
of $V$.

Alternatively, the weights may be computed by regarding $V$ as a
$\U_K$-module. Here the zero part $\U_K^0$ acts semisimply on $V$.
Moreover, $\U_K^0 = \U_\Z^0 \otimes_\Z K$ and $\U_\Z^0$ is generated
by all $\binom{H_i}{m}$ $(m\ge0)$, so the action of $\U_K^0$ on any
$V$ is determined by the action of $H_i \otimes 1$, for $i =
1,\dots,n$. So $V = \bigoplus_\lambda V_\lambda$ where $V_\lambda = \{
v\in V \mid (H_i \otimes 1)v = \lambda_i \,v, \text{ for all }
i=1,\dots,n \}$. One easily checks that the weights (and weight
spaces) in this sense coincide with the weights in the sense of the
preceding paragraph.

A weight $\lambda \in \Z^n$ is {\em dominant} (relative to the Borel
subgroup of upper triangular matrices in $\Gamma$) if $\lambda_1 \ge
\lambda_2 \ge \cdots \ge \lambda_n$. The Weyl group $\Sym_n$ acts on
$\Z^n$ by index permutation, {\em i.e.}\ $\sigma\,\lambda =
(\lambda_{\sigma^{-1}(1)}, \dots, \lambda_{\sigma^{-1}(n)})$, and the
set of dominant weights is a set of representatives for the orbits.

Not every element of $\Z^n$ can appear in the set of weights of an
object of $\M_K(n;r,s)$.  The following result describes those weights
that can appear.

\begin{lem} \label{lem:wts}
Let $\Lambda(n;r,s)$ be the set of all $\lambda\in \Z^n$ such that
$\sum \{\lambda_i: \lambda_i > 0\} = r-t$ and $\sum \{\lambda_i:
\lambda_i < 0\} = t-s$ for some $t$, $0 \le t \le \min(r,s)$. Let
$\Lambda^+(n;r,s)$ be the set of dominant weights in $\Lambda(n;r,s)$.

\par\noindent(a) Let $V$ be an object of $\M_K(n;r,s)$. The set of
weights of $V$ is contained in $\Lambda(n;r,s)$. The set of weights of
${\E_K}^{r,s}$ is $\Lambda(n;r,s)$.

\par\noindent(b) Let $\pi$ be the set of dominant weights of
${\E_K}^{r,s}$, as in the preceding section. Then $\pi =
\Lambda^+(n;r,s)$.
\end{lem}

\begin{proof}
Let $\varepsilon_1, \dots, \varepsilon_n$ be the standard basis of
$\Z^n$.  The weight of $\v_i$ is $\varepsilon_i$, so the weight of
$\v'_i$ is $-\varepsilon_i$. Thus, in the notation of \ref{Inrs}, for
a pair $(I,J)\in \I(n;r,s) = \I(n,r) \times \I(n,s)$ the weight of
\[
\v_{I,J}:= \v_{i_1}\otimes\cdots\otimes \v_{i_r}\otimes
\v'_{j_1}\otimes\cdots\otimes \v'_{j_s}
\]
is $\lambda \in \Z^n$ such that $\lambda_i$ is the difference between
the number of occurrences of $\v_i$ and the number of occurrences of
$\v'_i$ in the tensor $\v_{I,J}$. So such tensors form a basis of
weight vectors for ${\E_K}^{r,s}$.  The second claim in part (a) follows
immediately. And part (b) follows from the second claim in (a). 

Now let $V$ be an object of $\M_K(n;r,s)$. The weights of $V$ are
independent of characteristic, so without loss of generality we may
assume that $K=\C$.  Then $V$ is completely reducible (all rational
$\GL_n(\C)$-modules are completely reducible).

The simple $S_\C(n;r,s)$-modules are the simple factors of
${\E_\C}^{r,s}$, so their highest weights lie in $\pi =
\Lambda^+(n;r,s)$ and their weights lie in $\Lambda(n;r,s)$ since
$\Lambda(n;r,s)$ is stable under the Weyl group. Thus the weights of
$V$ lie in $\Lambda(n;r,s)$. This proves the first claim in (a).
\end{proof}

\begin{rmk}
The description of $\pi = \Lambda^+(n;r,s)$ in the preceding lemma can
be used to give a combinatorial proof that the set $\pi$ is saturated,
in the sense used by Donkin \cite{Donkin:SA1}. This means that if $\mu
\in \pi$ and $\lambda \unlhd \mu$ (in the dominance order on $\Z^n$;
see \S\ref{cellularbasis} ahead) for some dominant $\lambda$, then
$\lambda \in \pi$. That this property holds for the set $\pi$ is
obvious, however. Indeed, $\pi$ is the set of weights of the module
${\E_\C}^{r,s}$ in $\M_\C(n;r,s)$, and thus is a union of sets of
weights of its simple factors. That the set of weights of a finite
dimensional simple rational $\C\Gamma$-module is saturated is well
known, and unions of saturated sets are clearly saturated.
\end{rmk}

There is an alternative description of the set $\Lambda(n;r,s)$. Let
$\lambda$ be an element of $\Z^n$. A {\em proper partial sum} of
$\lambda_i$'s is a sum of the form $\lambda_{i_1} + \lambda_{i_2} +
\cdots + \lambda_{i_k}$, where $i_1, i_2, \dots, i_k$ are distinct
integers in the interval $[1,n]$, $0 < k < n$.

\begin{lem} \label{alt}
$\Lambda(n;r,s)$ is the set of all $\lambda \in \Z^n$ satisfying the
conditions
\par\noindent(a) $\textstyle\sum \lambda_i = r-s$; 
\par\noindent(b) $P \in [-s, r]$, for all proper partial sums $P$ of
$\lambda_i$'s.
\end{lem}

\begin{proof}
Suppose $\lambda\in \Lambda(n;r,s)$.  Then $\lambda \in \Z^n$ and
$\sum \{\lambda_i: \lambda_i>0\} = r-t$, $\sum \{\lambda_i:
\lambda_i<0\} = t-s$, for some $0 \le t \le \min(r,s)$.  In
particular, $\sum \lambda_i = r-s$.  Let $P$ be a proper partial sum
of $\lambda_i$'s. Write $P = P^+ + P^-$ where $P^+$ is the sum of the
positive contributions to the sum $P$ and $P^-$ is the sum of the
negative contributions.  If there are no positive (resp., negative)
summands in $P$, then $P^+$ (resp., $P^-$) is defined to be $0$. Then
one easily sees that $P^+ \in [0,r]$, $P^- \in [-s,0]$, and so $P \in
[-s, r]$.

On the other hand, if $\lambda$ satisfies conditions (a), (b) above,
then in particular $\sum \{\lambda_i: \lambda_i > 0\} \in [0,r]$ and
$\sum \{\lambda_i: \lambda_i < 0\} \in [-s,0]$. Consequently, $\sum
\{\lambda_i: \lambda_i > 0\} = r-t$ for some $0\le t \le r$ and $\sum
\{\lambda_i: \lambda_i < 0\} = t'-s$ for some $0 \le t' \le s$. But
putting the two sums together must give $r-s$, so $t=t'$ and $0 \le t
\le \min(r,s)$. The proof is complete.
\end{proof}

\subsection{}
There is a bijection between dominant weights and certain pairs of
partitions. Given a dominant weight $\lambda \in \Z^n$, let a pair
$(\lambda^+,\lambda^-)$ of partitions be determined as follows:
\begin{enumerate}
\item $\lambda^+ = (\lambda_1, \dots, \lambda_i)$ where $\lambda_i$ is
the rightmost positive entry of $\lambda$; if $\lambda$ has no
positive entries then $\lambda^+$ is the empty partition.
\item $\lambda^- = (-\lambda_n, \dots, -\lambda_j)$ where $\lambda_j$
is the leftmost negative entry of $\lambda$; if $\lambda$ has no
negative entries then $\lambda^-$ is the empty partition.
\end{enumerate}
In other words, $\lambda^+$ is the partition of the positive entries
in $\lambda$ and $\lambda^-$ is the partition obtained from the
negative entries by writing their absolute values in reverse order.
The reader may easily check that this procedure defines a bijection
between the set of dominant weights and the set of pairs of partitions
of total length not exceeding $n$. The element of $\Z^n$ corresponding
to a given pair of partitions (of total length not exceeding $n$) is
obtained by following the first partition with the negative reverse of
the second, inserting as many zeros in between as needed to make an
$n$-tuple.

Under this bijection, the set $\pi=\Lambda^+(n;r,s)$ corresponds with
the set of all pairs of partitions such that for some $0 \le t \le
\min(r,s)$ the first member of the pair is a partition of $r-t$, the
second is a partition of $s-t$, and the combined number of parts in
the pair does not exceed $n$.

\section{$S_K(n;r,s)$ is a generalized Schur algebra}\label{4}\noindent
In this section, we prove that the rational Schur algebra $S_K(n;r,s)$
may be identified with Donkin's generalized Schur algebra $S_K(\pi)$,
where $\pi=\Lambda^+(n;r,s)$ is the (saturated) set of dominant
weights occurring in ${\E_K}^{r,s} = {\E_K}^{\otimes r} \otimes
{\E_K^*}^{\otimes s}$.

\subsection{} \label{Zform}
We recall from \cite[3.2]{Donkin:SA1} that $S_R(\pi)$, for any
integral domain $R$, is constructed as $\U_R/I_R$ where $\U_R =
\U_\Z\otimes_\Z R$ is the hyperalgebra taken over $R$ and where $I_R$
is the ideal of $\U_R$ consisting of all elements of $\U_R$ that
annihilate every admissible $\U_R$-module belonging to $\pi$. A
$\U_R$-module is admissible if it is finitely generated and free over
$R$; such modules are direct sums of their weight spaces. (Weight
spaces are computed using the action of the zero part $\U_R^0=\U_\Z^0
\otimes_\Z R$ of $\U_R$, which is a set of commuting semisimple
operators.) A $\U_R$-module belongs to $\pi$ if every dominant weight
of the module lies in the set $\pi$.

In particular, taking $\pi = \Lambda^+(n;r,s)$, the set of dominant
weights of ${\E_\Q}^{r,s}$, we have $S_R(\pi)=\U_R/I_R$.  In case
$R=\Z$ we denote this algebra by $S_\Z(n;r,s)$. That it is an integral
form for the rational Schur algebra is the content of part (b) of the
next result.

\begin{thm} \label{thm:recon}
Fix $n\ge 2$, $r,s \ge 0$ and let $\pi = \Lambda^+(n;r,s)$ be the set
of dominant weights occurring in ${\E_\Q}^{r,s}$.  Let $K$ be an
infinite field.
\par\noindent(a) $S_K(n;r,s) = S_K(\pi)$.
\par\noindent(b) $S_K(n;r,s) \simeq S_\Z(n;r,s)\otimes_\Z K$. 
\end{thm}

\begin{proof}
(a) Let $I'_K$ be the kernel of the representation $\varphi_K: \U_K
\to \End_K({\E_K}^{r,s})$ (see \eqref{phi_K}). Since ${\E_K}^{r,s}$
belongs to $\pi$, $I_K$ is contained in $I'_K$.

Since $\E_K$, $\E_K^*$ are irreducible Weyl modules, they are tilting
modules.  Thus ${\E_K}^{r,s}$ is a tilting module. (The category of
tilting modules is closed under tensor products; see \cite[Part II,
E.7]{Jantzen}.) Hence ${\E_K}^{r,s}$ has a $\Delta$-filtration; that is, a
series of submodules
\[
  0 = F_0 \subset F_1 \subset \cdots \subset F_M = {\E_K}^{r,s}
\]
such that $F_{i+1}/F_i \simeq \Delta(\lambda_i)$ for each $i$. Here
$\Delta(\lambda_i)$ is the Weyl module of highest weight $\lambda_i
\in \pi$. 

In characteristic zero, every $\Delta(\lambda)$ for any $\lambda \in
\pi$ occurs as a filtration quotient (in fact, as a direct summand) of
a $\Delta$-filtration; see \cite{Stembridge}.  Thus for any $K$ and
any $\lambda \in \pi$ there must be some index $i$ such that $\lambda
= \lambda_i$; in other words, $\Delta(\lambda)$ is a sub-quotient of
${\E_K}^{r,s}$. It follows that the irreducible module $L(\lambda)$ of
highest weight $\lambda$ is a sub-quotient of ${\E_K}^{r,s}$, and that
any element of $I'_K$ annihilates $L(\lambda)$.

Hence any element of $I'_K$ must annihilate every irreducible
$\U_K$-module belonging to $\pi$. This shows the opposite inclusion
$I'_K \subset I_K$. Thus $I_K = I'_K$ and $S_K(n;r,s) = S_K(\pi)$.

(b) This follows from \cite[(3.2b)]{Donkin:SA1}.  Although in that
reference the algebraic group is semisimple, the results are equally
valid in the reductive case. Alternatively, one can use \ref{SL-cor}
to reduce the question to the semisimple case, and then apply
\cite[(3.2b)]{Donkin:SA1}.
\end{proof}

\subsection{}
Note that Donkin \cite[(3.2b)]{Donkin:SA1} showed that the natural map
$I_\Z \otimes_\Z R \to I_R$ induces an isomorphism
$S_\Z(\pi)\otimes_\Z R \to S_R(\pi)$, for any integral domain $R$.
Thus it makes sense to define $S_R(n;r,s)$ for any integral domain $R$
by $S_R(n;r,s) = S_\Z(n;r,s)\otimes_\Z R$. In particular, $S_K(n;r,s)$
is now defined for any field $K$, not just for infinite fields.

As an immediate consequence of the preceding theorem, we obtain the
result that rational Schur algebras are quasihereditary, in the sense
defined by Cline, Parshall, and Scott \cite{CPS}. In particular, this
means that rational Schur algebras have finite global dimension.

\begin{cor}
$S_K(n;r,s)$ is quasihereditary, for any field $K$.
\end{cor}

\begin{proof}
It easily follows from results of \cite{Donkin:SA1} that any
generalized Schur algebra (over a field) is quasihereditary.
\end{proof}

An alternate proof of the preceding result may be given based on the
cellular structure described in the next section.

\section{Cellular bases for $S_K(n;r,s)$}\label{5a}\noindent
The purpose of this section is to show that every rational Schur
algebra $S_K(n;r,s)$ inherits a cellular basis from a certain ordinary
Schur algebra, of which it is a quotient.

\subsection{Cellular algebras}

For the reader's convenience, we recall the original definition of
cellular algebra from Graham and Lehrer \cite{GL}. 

Let $A$ be an associative algebra over a commutative ring $R$, free
and of finite rank as an $R$-module. The algebra $A$ is {\em cellular}
with cell datum $(\pi,M,C,\iota)$ if the following conditions are
satisfied:

(C1) The finite set $\pi$ is a partially ordered set, and for each
$\lambda \in \pi$ there is a finite set $M(\lambda)$, such that
the algebra $A$ has an $R$-basis $\{ C^\lambda_{S,\,T} \}$, where
$(S,T)$ runs over all elements of $M(\lambda)\times M(\lambda)$ and
$\lambda$ runs over $\pi$.

(C2) The map $\iota$ is an $R$-linear anti-involution of $A$,
interchanging $C^\lambda_{S,\,T}$ and $C^\lambda_{T,\,S}$.

(C3) For each $\lambda \in \pi$, $a\in A$, and $S,T \in
M(\lambda)$ we have
\[
a C^\lambda_{S,\,T} = \sum_{U \in M(\lambda)}
r_a(U,S)C^\lambda_{U,\,T}       \qquad(r_a(U,S)\in R)
\]
modulo a linear combination of basis elements with upper index $\mu$
strictly less than $\lambda$ in the given partial order.

Note that the coefficients $r_a(U,S)$ in the expression in (C3) do not
depend on $T \in M(\lambda)$.

Suppose that $A$ is cellular. For each $\lambda \in \pi$, let $A[\le
\lambda]$ (respectively, $A[< \lambda]$) be the $R$-submodule of $A$
spanned by all $C^\mu_{S,\,T}$ such that $\mu \le \lambda$
(respectively, $\mu < \lambda$) and $S,T \in M(\mu)$. It is easy to
see that $A[\le \lambda]$ and $A[< \lambda]$ are two-sided ideals of
$A$, for each $\lambda \in \pi$. If we choose some total ordering for
the elements of $\pi$, say $\pi=\{\lambda^{(1)}, \dots,
\lambda^{(t)}\}$, refining the given partial order $\le$ on $\pi$,
setting $A_i = A[\le \lambda^{(i)}]$ gives a chain of two-sided ideals
of $A$
\begin{equation}
\{0\} \subset A_1 \subset A_2 \subset \cdots \subset A_t = A
\end{equation}
such that $\iota(A_i)=A_i$ for each $1\le i\le t$.

Axiom (C3) guarantees that, for any $T \in M(\lambda)$, the
$R$-submodule of $A[\le \lambda]/A[< \lambda]$ spanned by all
\begin{equation} \label{std:basis}
C^\lambda_{S,\,T} + A[<\lambda] \qquad (S\in M(\lambda)) 
\end{equation}
is a left $A$-module. We denote this $A$-module by $\Delta^\lambda_T$.
The elements in \eqref{std:basis} form an $R$-basis of this module, so
$\Delta^\lambda_T$ is free over $R$ of $R$-rank $|M(\lambda)|$. Thus
$A[\le \lambda]/A[< \lambda]$ is the direct sum of $\Delta^\lambda_T$
as $T$ runs over $M(\lambda)$; all these left $A$-modules are
isomorphic to the abstract left $A$-module $\Delta(\lambda)$ with
basis $C_S$ ($S\in M(\lambda)$) and action given by $aC_S =\sum_{U \in
M(\lambda)} r_a(U,S)C_U$ for any $a \in A$. The module
$\Delta(\lambda)$ is called a cell module.

By applying the involution $\iota$ we may decompose $A[\le
\lambda]/A[< \lambda]$ into a similar direct sum of $|M(\lambda)|$
right $A$-modules, all of which are isomorphic with
$\iota(\Delta(\lambda))$.  In fact, the $A$-bimodule $A[\le
\lambda]/A[< \lambda]$ is isomorphic with $\Delta(\lambda) \otimes_R
\iota(\Delta(\lambda))$. Note that \cite{KX} give an alternative
approach to the theory of cellular algebras based on these bimodules.

\subsection{Cellular bases of $S_K(n,r)$} \label{cellularbasis}
There are two well known cellular bases of a classical Schur algebra
$S_K(n,r)$: the codeterminant basis of J.A.~Green \cite{Green:codet}
and the canonical basis (see \cite{BLM}; \cite{Du}).  The former is
compatible with the Murphy basis \cite{Murphy:1, Murphy:2} and the
latter with the Kazhdan--Lusztig basis \cite{KL} (for the symmetric
group algebra or the Hecke algebra in type $A$). Both bases of
$S_K(n,r)$ may be obtained by lifting the corresponding basis from the
group algebra of the symmetric group. Alternatively, the latter basis
may be realized by a descent from Lusztig's modified form (see
\cite{Lusztig:book}) of $\U_\Q$. (See \cite[\S6.14]{Doty} for details.)

In both cases, one takes the set $\pi$ to be $\Lambda^+(n,r)$;
this may be identified with the set of partitions of $r$ into not more
than $n$ parts. The partial order on $\pi$ is the {\em reverse}
dominance order; {\em i.e.}, we have to read $\le$ as $\unrhd$. Here
$\mu \unrhd \lambda$ (for any $\mu = (\mu_1, \dots, \mu_n)$, $\lambda
= (\lambda_1, \dots, \lambda_n) \in \Z^n$) if $\sum_{1\le j\le i}
\mu_j \ge \sum_{1\le j\le i} \lambda_j$ for all $i=1, \dots, n$.  This
is the partial order (with respect to the upper triangular Borel
subgroup) on the set of weights $\Z^n$ defined by $\mu \unrhd \lambda$
if $\mu - \lambda$ is a sum of positive roots.  For each $\lambda \in
\pi$, the set $M(\lambda)$ is the set of row semistandard
$\lambda$-tableaux, {\em i.e.}, Young diagrams of shape $\lambda$ with
entries from the set $\{1, \dots, n \}$ such that entries strictly
increase down columns and weakly increase along rows. In both cases
there is a cellular basis $\{ C^\lambda_{S,\,T} \}$ of $S_K(n,r)$, and
in both cases the cell modules $\Delta(\lambda)$ are isomorphic with
the Weyl modules of highest weight $\lambda \in \pi$.

The codeterminant basis of $S_K(n,r)$ is quite simple to describe.
Let $\xi^{}_{I,\,J}$ for $I,J \in \I(n,r)$ be the elements of
$S_K(n,r)$ described in \cite[\S2.3]{Green:book}.  For each pair $S,
T$ of row semistandard $\lambda$-tableaux, for $\lambda \in
\Lambda^+(n,r)$, one obtains corresponding elements $I,J \in \I(n,r)$
by reading the entries in left-to-right order across the rows of each
tableau, in order from the top row to the bottom row. Then
$C^{\lambda}_{S,\,T} = \xi^{}_{I,\,\ell} \,\xi^{}_{\ell,\,J}$ where
$\ell=\ell(\lambda)$ is the element of $\I(n,r)$ corresponding to the
$\lambda$-tableau with all entries in the $i$th row equal to $i$, as
$i$ varies through the rows. That this is in fact a cellular basis
follows easily from the straightening algorithm of Woodcock
\cite{Wo}. A $q$-analogue of this straightening algorithm was given in
\cite{RMG:straight} and a proof of the cellularity of the
codeterminant basis can be found in \cite[Proposition
6.2.1]{RMG:thesis}; that same argument is valid in the $q=1$ case.

\subsection{The quotient map $S_K(n,r+(n-1)s) \to S_K(n;r,s)$}

Let $\mathcal{C}$ be the category of rational $K\Gamma$-modules.
There is a functor $\Psi$ from $\mathcal{C}$ to $\mathcal{C}$, sending
$V$ to $V\otimes \det$. This is clearly an invertible functor. For any
$s\in \Z$ we have a corresponding functor $\Psi^s$ sending $V$ to
$V\otimes \det^{\otimes s}$ if $s\ge 0$, and sending $V$ to $V \otimes
(\det^{-1})^{\otimes |s|}$ in case $s<0$. Clearly $\Psi^{-s}$ is
inverse to $\Psi^s$ for any $s\in \Z$. In fact, $\Psi^r \circ \Psi^s =
\Psi^{r+s}$ for any $r,s \in \Z$, and $\Psi^0$ is the identity
functor.

Set $X = \Z^n$, regarded as additive abelian group, and consider its
group ring $\Z[X]$, with basis $\{ e(\lambda) \mid \lambda \in X\}$
and product determined by $e(\lambda)\cdot e(\mu) = e(\lambda+\mu)$
for $\lambda,\mu \in X$. The formal character of $V$ is the element of
$\Z[X]$ given by $\ch V = \sum_{\lambda \in X} (\dim
V_\lambda)e(\lambda)$, where $V_\lambda$ is the $\lambda$-weight space
as in \S\ref{wts}.  Clearly we have
\begin{equation}
  \ch (\Psi^s V) = \sum_{\lambda\in X} (\dim
  V_\lambda)\,e(\lambda+s\omega)
\end{equation}
for any $V$ and any $s$.  Here $\omega:=(1^n)=(1,\dots,1) \in \Z^n$.
In particular, if $V$ is a highest weight module of highest weight
$\lambda$ then $\Psi^s V$ is a highest weight module of highest weight
$\lambda+s\omega$ (for any $s\in \Z$).

The existence of the quotient map is motivated by the observation that
$\Psi \E^*_K = \E^*_K \otimes \det \simeq \Lambda^{n-1} \E_K$ (as
$K\Gamma$-modules or $\U_K$-modules). The existence of this
isomorphism is clear from comparing the highest weight of each module
(both are irreducible in any characteristic).  Thus we have an
isomorphism of rational $K\Gamma$-modules
\begin{equation}
\Psi^s({\E_K}^{\otimes r} \otimes {\E_K^*}^{\otimes s}) \simeq
 {\E_K}^{\otimes r} \otimes (\Lambda^{n-1} \E_K)^{\otimes s}
\end{equation}
for any $r,s \ge 0$.  The $K\Gamma$-module $\Lambda^{n-1} \E_K$ is the
Weyl module of highest weight $(1,\dots,1,0) \in \Z^n$, and may be
realized as a submodule of tensor space ${\E_K}^{\otimes (n-1)}$ by
means of the Carter--Lusztig construction, as the submodule spanned by
all anti-symmetric tensors of the form
\[
  \sum_{\sigma \in \Sym_{n-1}} (-1)^{\sgn(\sigma)} \v_{\sigma(i_1)}
  \otimes \cdots \otimes \v_{\sigma(i_{n-1})}
\]
for any $I = (i_1, \dots, i_{n-1}) \in \I(n,n-1)$. See \cite[5.2,
Example 2]{Green:book} for details. 

The above embedding $\gamma:\Lambda^{n-1} \E_K \to {\E_K}^{\otimes
(n-1)}$ induces a corresponding embedding
\begin{equation}
  {\E_K}^{\otimes r} \otimes (\Lambda^{n-1} \E_K)^{\otimes s} 
  \to {\E_K}^{\otimes r} \otimes {\E_K}^{\otimes (n-1)s} = 
  {\E_K}^{\otimes (r+(n-1)s)}
\end{equation}
obtained by using the identity map in the first $r$ tensor positions
and repeating $\gamma$ in the rest of the factors $s$ times. 

By results in \S\ref{2} we have an injective map $\AA_K(n;r,s)
\to A_K(n,r+(n-1)s)$ given by $c \to c\cdot \d^s$ for any $c \in
\AA_K(n;r,s)$. Since $\d$ is a group-like element, this map is a
coalgebra morphism. By dualizing we therefore obtain a surjective
algebra map
\begin{equation} \label{quomap}
  S_K(r+(n-1)s) \to S_K(n;r,s)
\end{equation}
for any $r,s\ge 0$ and any field $K$. 

When $n=2$ the rational Schur algebras are not new. More precisely,
we have the following result.

\begin{prop}\label{n=2case}
The quotient map $S_K(2,r+s) \to S_K(2;r,s)$ is an isomorphism of
algebras, for any $r,s\ge 0$.
\end{prop}

\begin{proof}
This follows from the isomorphism $\E^*_K \otimes \det \simeq \E_K$
(when $n=2$). This implies that $\Psi^s {\E_K}^{r,s} \simeq
{\E_K}^{\otimes (r+s)}$, and it follows by taking coefficient spaces
that $\AA_K(2;r,s) \to A_K(2,r+s)$ is an isomorphism of
coalgebras. The result follows by dualizing.
\end{proof}

Return now to general $n \ge 2$, and let $r,s \ge 0$ be given.  Given
any $S_K(n;r,s)$-module $V$, we may regard $V$ as an
$S_K(n,r+(n-1)s)$-module by composing the action with the quotient map
$S_K(n,r+(n-1)s) \to S_K(n;r,s)$.  Denote by $V'$ the resulting
$S_K(n,r+(n-1)s)$-module.  It is easy to see that $V'$ is isomorphic
with $\Psi^s V = V\otimes \det^{\otimes s}$, as
$S_K(n,r+(n-1)s)$-modules. Since the weights
$\lambda=(\lambda_1,\dots,\lambda_n) \in \Z^n$ of $V$ satisfy the
condition $\lambda_i \in [-s,r]$ for each $i = 1,\dots, n$ it is clear
that the weights $\mu=(\mu_1,\dots,\mu_n) \in \Z^n$ of $\Psi^s V$
satisfy the condition $\mu_i \in [0,r+s]$ for each $i=1,\dots,n$.

Note that if $V$ is a highest weight module for $S_K(n,r+(n-1)s)$ of
highest weight $\lambda$, then $\Psi^{-s} V$ is a highest weight
module for $S_K(n;r,s)$ if and only if the weight $\lambda$ satisfies
the condition $\lambda_i \in [0,r+s]$ for each $i = 1,\dots, n$.

\begin{thm}\label{thm:descent} 
Let $\{ C^\lambda_{S,\,T} \}$ be any cellular basis of the ordinary
Schur algebra $S_K(n, r+(n-1)s)$ such that $\pi=\Lambda^+(n,
r+(n-1)s)$ ordered by the reverse dominance order $\unrhd$.  Then the
kernel of the quotient map $S_K(n,r+(n-1)s) \to S_K(n;r,s)$ is spanned
by the set of all $C^\mu_{S,\,T}$ such that $\mu$ has at least one
part exceeding $r+s$.
\end{thm}

\begin{proof}
Set $A = S_K(n,r+(n-1)s)$.  Set $J$ equal to the sum of all ideals
$A[\unrhd \mu]$ for $\mu \in \pi$ satisfying $\mu_i > r+s$ for some
index $i$. The set of such $\mu$ is a saturated subset $\pi'$ of
$\pi$ under the dominance order, and the two-sided ideal $J$,
regarded as a left $A$-module, is filtered by $\Delta(\nu)$ for $\nu
\in \pi'$. Thus $J$ is contained in the kernel, so $S_K(n;r,s)$ is
a quotient of $A/J$.

The quotient $A/J$ has a filtration by $\Delta(\lambda)$ with $\lambda
\in \pi-\pi'$, and its dimension is $\sum_\lambda (\dim
\Delta(\lambda))^2$ as $\lambda$ runs over the set
$\pi-\pi'$. But $\dim \Delta(\lambda) = \dim \Psi^{-s}
\Delta(\lambda) = \Delta(\lambda-s\omega)$, so $\dim A/J = \sum_\mu
(\dim \Delta(\mu))^2$ as $\mu$ runs over the set $\Lambda^+(n;r,s)$.
Thus $A/J$ has the same dimension as $S_K(n;r,s)$. The proof is
complete.
\end{proof}

\begin{cor}
Under the quotient map $S_K(n,r+(n-1)s) \to S_K(n;r,s)$, the images of
all $C^\lambda_{S,\,T}$, for $\lambda$ satisfying the condition
$\lambda_i\le r+s$ (for $i = 1, \dots, n$) form a cellular basis of
the rational Schur algebra $S_K(n;r,s)$.
\end{cor}

To multiply two basis elements obtained by the procedure in the
corollary, one computes their product in $A=S_K(n,r+(n-1)s)$ expressed
as a linear combination of the other cellular basis elements, omitting
any terms belonging to the ideal $J = \sum_{\mu \in \pi'} A[\unrhd
\mu]$ (in the notation of the proof of \ref{thm:descent}).

Note that if the cellular basis of $S_K(n,r+(n-1)s)$ is defined over
$\Z$ then the same is true of the basis inherited by the quotient
$S_K(n;r,s)$.  Since the codeterminant and the canonical bases of
Schur algebras both have this property, we see that in fact the above
procedure gives cellular bases for $S_\Z(n;r,s)$ which induce via
change of base ring the corresponding cellular basis of $S_K(n;r,s)$
for each field $K$.

\subsection{Examples}
As previously mentioned, $S_K(n;r,0)=S_K(n,r)$ for any $r$.  By
\ref{opp} we have an isomorphism $S_K(n;0,s) \simeq S_K(n,s)^{\mathrm{opp}}$
for any $s$.  Moreover, by \ref{n=2case} we have $S_K(2;r,s) \simeq
S_K(2,r+s)$ for any $r,s$. 

\smallskip

Thus the smallest interesting new example is $S_K(3;1,1)$. According
to \ref{thm:descent}, this is a quotient of $S_K(3,3)$. The Weyl
modules for $S_K(3,3)$ are 
\[
\begin{array}{ll}
\Delta(3,0,0) & \dim=10\\
\Delta(2,1,0) & \dim=8\\
\Delta(1,1,1) & \dim=1.
\end{array}
\]
(The dimension of $\Delta(\lambda)$, for a partition $\lambda \in
\Lambda^+(n,r)$, is the number of semistandard $\lambda$-tableaux.)
So $S_K(3,3)$ has dimension $165=10^2+8^2+1^2$. The quotient map
$A=S_K(3,3) \to S_K(3;1,1)$ has kernel $A[\unrhd (3,0,0)]$; the kernel
has dimension 100 and is spanned by all $C^\lambda_{S,\,T}$ for
$\lambda=(3,0,0)$ as $S,T$ varies over all pairs of
$\lambda$-tableaux. Thus $S_K(3;1,1)$ has dimension 65, and has two
Weyl modules
\[
\begin{array}{ll}
\Delta(1,0,-1) = \Psi^{-1}\Delta(2,1,0) & \dim=8\\
\Delta(0,0,0) = \Psi^{-1}\Delta(1,1,1) & \dim=1.
\end{array}
\]
One could at this point write out the elements of a cellular basis for
$S_K(3;1,1)$, indexed by pairs of semistandard $\mu$-tableaux for $\mu
= (2,1,0)$ and $(1,1,1)$, and compute the structure constants with
respect to the basis.

\smallskip

The next simplest cases are $S_K(3;2,1)$ or $S_K(3;1,2)$. These are
related by $S_K(3;2,1) \simeq S_K(3;1,2)^{\mathrm{opp}}$; both algebras have
dimension 270.  Note that $S_K(3,2,1)$ is a quotient of $S_K(3,4)$
while $S_K(3;1,2)$ is a quotient of $S_K(3,5)$, according to
\ref{thm:descent}. 

Now $\Lambda^+(3,4) = \{(4,0,0), (3,1,0), (2,2,0), (2,1,1)\}$. The
Weyl modules for $S_K(3,4)$ are
\[
\begin{array}{ll}
  \Delta(4,0,0) & \dim=15\\
  \Delta(3,1,0) & \dim=15\\
  \Delta(2,2,0) & \dim=6\\
  \Delta(2,1,1) & \dim=3.
\end{array}
\]
The dimension of $S_K(3,4)$ is $495 = 15^2 + 15^2 + 6^2 + 3^2$.  The
kernel of the quotient map $A=S_K(3,4) \to S_K(3;2,1)$ is $A[\unrhd
(4,0,0)]$ of dimension 225. Thus the Weyl modules for $S_K(3;2,1)$ are
\[
\begin{array}{ll}
  \Delta(2,0,-1)= \Psi^{-1}\Delta(3,1,0) & \dim=15\\
  \Delta(1,1,-1)= \Psi^{-1}\Delta(2,2,0) & \dim=6\\
  \Delta(1,0,0) = \Psi^{-1}\Delta(2,1,1) & \dim=3.
\end{array}
\]
Again, one could write out a cellular basis at this point.

Similarly, $\Lambda^+(3,5) = \{(5,0,0), (4,1,0), (3,2,0), (3,1,1),
(2,2,1)\}$ and the Weyl modules for $S_K(3,5)$ are
\[
\begin{array}{ll}
\Delta(5,0,0) & \dim=21\\
\Delta(4,1,0) & \dim=24\\
\Delta(3,2,0) & \dim=15\\
\Delta(3,1,1) & \dim=6\\
\Delta(2,2,1) & \dim=3.
\end{array}
\]
In this case $\dim S_K(3,5)$ is $1287 = 21^2 + 24^2 + 15^2 + 6^2 +
3^2$ and the kernel of the quotient map $A=S_K(3,5) \to S_K(3;1,2)$ is
$A[\unrhd (4,1,0)]$ of dimension $1017 = 21^2+24^2$. Note that the kernel
contains $A[\unrhd (5,0,0)]$.  The Weyl modules for $S_K(3;1,2)$ are
\[
\begin{array}{ll}
\Delta(1,0,-2) = \Psi^{-2}\Delta(3,2,0) & \dim=15\\
\Delta(1,-1,-1) = \Psi^{-2}\Delta(3,1,1)& \dim=6\\
\Delta(0,0,-1) =\Psi^{-2}\Delta(2,2,1) & \dim=3.
\end{array}
\]
As above, one could write out a cellular basis at this point.

\smallskip

The next example would be $S_K(3;2,2)$, which has dimension 994 and is
realized as a quotient of $S_K(3,6)$ (of dimension 3003). We leave the
details to the reader.

\section{Generators and relations}\label{7}\noindent
Another consequence of the existence of the quotient map
\eqref{quomap} is a simple description of $S_\Q(n;r,s)$ by generators
and relations, generalizing the presentation of $S_\Q(n,d)$ obtained
in \cite{DG:PSA}. We recall the elements $\varepsilon_1, \dots,
\varepsilon_n$ from the proof of Lemma \ref{lem:wts}; these elements
form a natural basis for the set $\Z^n$ of weights of $\Gamma =
\GL_n(K)$. We set $\alpha_j = \varepsilon_j - \varepsilon_{j+1}$ for
$1 \le j \le n-1$. These are the usual simple roots in type
$A_{n-1}$. We have an inner product on $\Z^n$ given by on generators
by $(\varepsilon_i, \varepsilon_j) = \delta_{ij}$.

In the following we work over $\Q$, but in fact $\Q$ may be replaced
by any field of characteristic zero. 

\subsection{} \label{PSA:rels}
We recall from \cite[Theorem 1.1]{DG:PSA} that the algebra $S_\Q(n,d)$
is isomorphic with the associative algebra with 1 generated by symbols
$e_i$, $f_i$ ($1 \le i \le n-1$) and $H_i$ ($1 \le i \le n$) subject
to the relations

(a) $H_i H_j = H_j H_i$
 
(b) $e_i f_j - f_j e_i = \delta_{ij}(H_j - H_{j+1})$

(c) $H_i e_j - e_j H_i = (\varepsilon_i, \alpha_j) e_j$, \ \ 
$H_i f_j - f_j H_i = -(\varepsilon_i, \alpha_j) f_j$ 

(d) $e_i^2 e_j - 2 e_i e_j e_i + e_j e_i^2 = 0$, \ \  
$f_i^2 f_j - 2 f_i f_j f_i + f_j f_i^2 = 0$
\quad $(|i-j|=1)$

(e) $e_i e_j = e_j e_i$,\ \ 
$f_i f_j = f_j f_i$ \quad $(|i-j| \ne 1)$ 

(f) $H_1 + H_2 + \cdots + H_n  = d$

(g) $H_i(H_i-1)\cdots(H_i-d+1)(H_i-d) = 0$.

\noindent
We also recall that Serre proved that the enveloping algebra
$\U=\U_\Q(\gl_n)$ is the $\Q$-algebra on the same generators but
subject only to relations (a)--(e).  Thus the algebra $S_\Q(n,d)$ is a
homomorphic image of $\U$, by the map given on generators by $e_i \to
e_i$, $f_i \to f_i$, $H_i \to H_i$.

\begin{prop}
Assume that $n \ge 2$ and set $A = S_\Q(n,d)$ with $d = r+(n-1)s$. Let
$A[\pi'] = \cup_{\mu \in \pi'} A[\unrhd \mu]$ where $\pi'$ is the set
of all $\mu$ in $\Lambda^+(n,d)$ such that at least one part $\mu_i$
exceeds $r+s$. Then the quotient algebra $A/A[\pi']$ is isomorphic
with the algebra given by the same generators as in \S\ref{PSA:rels}
subject to relations \ref{PSA:rels}(a)--(f) along with 

(g') $H_i(H_i-1)\cdots(H_i-r-s+1)(H_i-r-s) = 0$.
\end{prop}

\begin{proof}
Since relation \ref{PSA:rels}(g) is a consequence of relation (g')
above it follows that the algebra $A'$ given by the generators and
relations of the proposition is a quotient of $A$. We let $I$ be the
kernel so that we have $A' \simeq A/I$.  Clearly $I$ is generated by
all $H_i(H_i-1)\cdots(H_i-r-s+1)(H_i-r-s)$.  To prove the proposition we
must show that $I$ coincides with the cell ideal $A[\pi']$.

To see this we first observe that the idempotents $1_\mu \in A$ (for
$\mu \in \pi'$) all map to zero in the quotient $A'$. These
idempotents were defined in \cite{DG:PSA} by $1_\mu = \prod_i
\binom{H_i}{\mu_i}$, but one easily sees as in \cite[Lemma 5.3]{DEH}
that they coincide with Green's idempotents $\xi_\mu$ defined in
\cite[\S3.2]{Green:book}. These are in fact codeterminants: $\xi_\mu =
\xi_{I, \ell} \xi_{\ell, I}$ where $I$ is an appropriate element of
$\I(n,r)$ ($\ell$ is the particular element of $\I(n,r)$ defined in
\ref{cellularbasis}). This is an easy consequence of the
multiplication rule \cite[(2.3b)]{Green:book}. Let $\xi = \sum
\xi_\mu$ where the sum is taken over all $\mu \in \pi'$. Then $\xi C =
C$ for all codeterminants $C$ in $A[\pi']$, so it follows that $C$
maps to zero in the quotient $A'$. Thus $A[\pi'] \subseteq I$.

On the other hand, every other idempotent maps to a nonzero idempotent
in the quotient $A'$, and thus it follows that every cell module
$\Delta(\lambda)$ indexed by some $\lambda \in \Lambda^+(n,d) - \pi'$
is an irreducible $A'$-module. Since $A'$ is a semisimple algebra it
follows that 
\[
\dim_\Q A' \ge \sum_{\lambda \in \Lambda^+(n,d) - \pi'}
(\dim_\Q \Delta(\lambda))^2.
\]
 But $\dim_\Q A' =
\dim_\Q A - \dim_\Q I$ and $\dim_\Q A = \sum_{\lambda \in
\Lambda^+(n,d)} (\dim_\Q \Delta(\lambda))^2$, so we get 
\[
\dim_\Q I \le \sum_{\mu \in \pi'} (\dim_\Q \Delta(\mu))^2 = \dim_\Q
A[\pi'].
\]  
Since $A[\pi'] \subseteq I$ we conclude that $I = A[\pi']$ by
dimension comparison. The proof is complete.
\end{proof}

\subsection{} \label{Snrs:rels}
The proof of the preceding proposition reveals that the irreducible
representations of $A' = A/A[\pi']$ are precisely the
$\Delta(\lambda)$ with $\lambda \in \Lambda^+(n,d)$ such that each
part $\lambda_i \in [0,r+s]$. As already noted in \S\ref{5a}, this set
of modules is in bijective correspondence with the set of irreducible
$S_\Q(n;r,s)$-modules.  Now we consider the effect of replacing the
generators $H_i$ in $A'$ with new generators $\tH_i = H_i - s$. This
induces an algebra automorphism of $A'$, and has the effect of
changing the defining relations slightly, to the following:

(a) $\tH_i \tH_j = \tH_j \tH_i$
 
(b) $e_i f_j - f_j e_i = \delta_{ij}(\tH_j - \tH_{j+1})$

(c) $\tH_i e_j - e_j \tH_i = (\varepsilon_i, \alpha_j) e_j$, \ \ 
$\tH_i f_j - f_j \tH_i = -(\varepsilon_i, \alpha_j) f_j$ 

(d) $e_i^2 e_j - 2 e_i e_j e_i + e_j e_i^2 = 0$, \ \  
$f_i^2 f_j - 2 f_i f_j f_i + f_j f_i^2 = 0$
\quad $(|i-j|=1)$

(e) $e_i e_j = e_j e_i$,\ \ 
$f_i f_j = f_j f_i$ \quad $(|i-j| \ne 1)$ 

(f) $\tH_1 + \tH_2 + \cdots + \tH_n  = r-s$

(g') $(\tH_i+s)(\tH_i+s-1)\cdots(\tH_i-r+1)(\tH_i-r) = 0$.

\noindent
This algebra should be regarded as a quotient of $\U$ via the map
given on generators by $e_i \to e_i$, $f_i \to f_i$, $H_i \to \tH_i +
s$. By the remarks preceding Theorem \ref{thm:descent}, the
irreducible modules for $A'$ now become precisely the set of
irreducible $\U$-modules whose highest weights belong to the set
$\Lambda^+(n;r,s)$. This gives the following result.

\begin{thm}
  For $n\ge 2$, the algebra $S_\Q(n;r,s)$ is isomorphic with the
  associative algebra with 1 given by the generators $e_i$, $f_i$
  $(1\le i \le n-1)$ and $\tH_i$ $(1\le i \le n)$ subject to the
  relations \ref{Snrs:rels}(a)--(g').
\end{thm}

Of course, this algebra is isomorphic with $A'$. By the theory of
semisimple algebras, it makes no difference whether we regard the
irreducible representations for the algebra as the $\Delta(\lambda)$
for $\lambda \in \pi'$ or the $\Delta(\lambda)$ for $\lambda \in
\Lambda^+(n;r,s)$.

\begin{rmk}
  In light of results of \cite{DG:PSA} it is now clear how one might
  define a natural candidate for the rational $q$-Schur algebra. It is
  the $\Q(v)$-algebra ($v$ an indeterminate) given by generators
  $E_i$, $F_i$ $(1\le i \le n-1)$ and $\tK_i$ $(1\le i \le n)$ with
  relations

(a) $\tK_i \tK_j = \tK_j \tK_i$
 
(b) $E_i F_j - F_j E_i = \delta_{ij} \frac{\tK_j {\tK_{j+1}}^{-1} - 
  {\tK_j}^{-1} \tK_{j+1} }{v - v^{-1}}$

(c) $\tK_i E_j - E_j \tK_i = v^{(\varepsilon_i, \alpha_j)} E_j$, \ \ 
$\tK_i F_j - F_j \tK_i = v^{-(\varepsilon_i, \alpha_j)} F_j$ 

(d) $E_i^2 E_j - (v+v^{-1}) E_i E_j E_i + E_j E_i^2 = 0$, \ \  
$F_i^2 F_j - (v+v^{-1}) F_i F_j F_i + F_j F_i^2 = 0$
\quad $(|i-j|=1)$

(e) $E_i E_j = E_j E_i$,\ \ 
$F_i F_j = F_j F_i$ \quad $(|i-j| \ne 1)$ 

(f) $\tK_1  \tK_2  \cdots  \tK_n  = v^{r-s}$

(g') $(\tK_i-v^{-s})(\tK_i-v^{-s+1})\cdots(\tK_i-v^{r-1})(\tK_i-v^r) = 0$. 
\end{rmk}

\section{Schur--Weyl duality}\label{6}\noindent
Let $K$ be a field of characteristic zero. Classical Schur-Weyl
duality, based on tensor space ${\E_K}^{r,0} = {\E_K}^{\otimes r}$,
has been extended in \cite{Benkart-et-al} to mixed tensor space
${\E_K}^{r,s} = {\E_K}^{\otimes r} \otimes {\E^*_K}^{\otimes s}$. Thus
we may regard the rational Schur algebra as a centralizer algebra for
the action of a certain algebra $\Br_{r,s}^{(n)}$, the so-called {\em
walled} or {\em rational} Brauer algebra.  This diagram algebra is
cellular \cite{GM}.

\subsection{The Brauer algebra $\Br_r^{(x)}$} \label{Brauer}

This algebra was introduced in Brauer \cite{Brauer}. 
Let $R$ be a commutative ring, let $R[x]$ be the ring of polynomials
in an indeterminate $x$, and consider the free $R[x]$-module
$\Br_r^{(x)}$ with basis consisting of the set of $r$-diagrams.  An
$r$-diagram is an (undirected) graph on $2r$ vertices and $r$ edges
such that each vertex is incident to precisely one edge.  One usually
thinks of the vertices as arranged in two rows of $r$ each, which are
then called the {top} and {bottom} rows.  Edges connecting two
vertices in the same row (different rows) are called {\em horizontal}
(resp., {\em vertical}).  We can compose two such diagrams $D_1$,
$D_2$ by identifying the bottom row of vertices in the first diagram
with the top row of vertices in the second diagram. The result is a
graph with a certain number, $\delta(D_1,D_2)$, of interior
loops. After removing the interior loops and the identified vertices,
retaining the edges and remaining vertices, we obtain a new
$r$-diagram $D_1 \circ D_2$, the {\em composite} diagram.
Multiplication of $r$-diagrams is then defined by the rule
\begin{equation}
D_1 \cdot D_2 = x^{\delta(D_1,D_2)} (D_1 \circ D_2).
\end{equation}
One can check that this multiplication makes $\Br_r^{(x)}$ into an
associative algebra.  Note that the subalgebra of $\Br_r^{(x)}$ spanned
by the diagrams containing just vertical edges may be identified with
the group algebra $R[x]\Sym_r$.

\subsection{The walled Brauer algebra}
The walled Brauer algebra (also called the rational Brauer algebra) is
a certain subalgebra of $\Br_{r+s}^{(x)}$ which first appeared in
\cite{Benkart-et-al}. Some examples and further details may be found
in \cite{Doty:SWD}.  

By definition, an $(r,s)$-diagram is an $(r+s)$-diagram in which we
imagine a wall separating the first $r$ columns of vertices from the
last $s$ columns of vertices, such that:

 (a) all horizontal edges cross the wall;

 (b) no vertical edges cross the wall.

\noindent
Let $\Br_{r,s}^{(x)}$ be the subalgebra of $\Br_{r+s}^{(x)}$ spanned
by the set of $(r,s)$-diagrams. This is a subalgebra of
$\Br_{r+s}^{(x)}$, with multiplication in $\Br_{r,s}^{(x)}$ the same
as that of $\Br_{r+s}^{(x)}$.  Given an $(r,s)$-diagram, interchanging
the vertices in each column on one side of the wall results in a
diagram with no horizontal edges.  Thus there is a bijection between
the set of $(r,s)$-diagrams and $\Sym_{r+s}$. It follows that
$\Br_{r,s}^{(x)}$ is a free $R$-module of rank $(r+s)!$.

Label the vertices on the top and bottom rows of an $(r,s)$-diagram by
the numbers $1, \dots, r$ to the left of the wall and $-1, \dots, -s$
to the right of the wall. Let $c_{i,-j}$ ($1 \le i \le r$; $1\le j \le
s$) be the diagram with a horizontal edge connecting vertices $i$ and
$-j$ on the top row and the same on the bottom row, and with all other
edges connecting vertex $k$ ($k \ne i, -j$) in the top and bottom
rows.  Given $\sigma \in \Sym_r$ let $t_\sigma$ denote the
$(r,s)$-diagram with the diagram corresponding to $\sigma$ appearing
to the left of the wall and the identity to the right.  Similarly,
given $\tau \in \Sym_s$ let $t'_\tau$ denote the $(r,s)$-diagram with
the diagram corresponding to $\tau$ appearing to the right of the wall
and the identity to the left.  It is easy to see that $\Br_{r,s}^{(x)}$
is generated by the permutations $\{t_\sigma, t'_\tau\}$ along with
just one of the $c_{i,-j}$.

\subsection{}
Now we specialize to $R=K$, a field, and set $x=n\cdot 1_K \in K$,
where as usual $n = \dim_K \E_K$.  The resulting $K$-algebra
$\Br_{r,s}^{(n)}$ acts on ${\E_K}^{r,s}$, on the right, as follows.  A
permutation diagram $t_\sigma$ $(\sigma \in \Sym_r)$ acts by place
permutation on the first $r$ factors of ${\E_K}^{r,s} = {\E_K}^{\otimes r}
\otimes {\E_K^*}^{\otimes s}$. Similarly a permutation diagram
$t'_\tau$ $(\tau \in \Sym_s)$ acts by place permutation on the last
$s$ factors of ${\E_K}^{r,s}$.  Finally, the element $c_{r,-1}$ acts as
a contraction in tensor positions $(r,-1)$; {\em i.e.},
\begin{equation}
\begin{aligned}
  (\v_{i_1}\otimes\cdots \otimes \v_{i_r} \otimes
  \v'_{j_1}\otimes\cdots \otimes \v'_{j_s})\, c_{r,-1} &\\ =
  \textstyle  
  \delta_{i_r,\,j_1}  \v_{i_1}\otimes\cdots
  \otimes\v_{i_{r-1}}\otimes (\sum_{k=1}^n \v_{k} &\otimes
  \v'_{k}) \otimes\v'_{j_2}\otimes \cdots \otimes \v'_{j_s}.
\end{aligned}
\end{equation}
The other $c_{i,-j}$ similarly contract in the $(i,-j)$th tensor
positions.

One checks that the above action of $\Br_{r,s}^{(n)}$ commutes with the
natural action of $\Gamma=\GL_n(K)$. 

\subsection{Double centralizer property} \label{DCP}
Assuming that $K$ is a field of characteristic zero and $n \ge r+s$,
it has been shown in \cite{Benkart-et-al} (see also \cite{Koike}) that
the two commuting actions on ${\E_K}^{r,s}$ mutually centralize one
another. In other words,
\begin{equation}\label{SWD}
\begin{gathered}
  \rho_K(K\Gamma ) = \End_{\Br_{r,s}^{(n)}}({\E_K}^{r,s}),\\
  \tau_K( \Br_{r,s}^{(n)} ) = \End_{\Gamma} ({\E_K}^{r,s})
\end{gathered}
\end{equation}
where $\tau_K$ is the representation of $\Br_{r,s}^{(n)}$ affording
the previously defined action of $\Br_{r,s}^{(n)}$ on ${\E_K}^{r,s}$.

When $K$ is any infinite field, we have the equality $S_K(n; r,s)
=\rho_K(K\Gamma)$, by \ref{rho-rep}.  It follows that the dimension of
$\rho_K(K\Gamma)$ is independent of the (infinite) field $K$.
Moreover, still for $K$ infinite, we also know that the dimension of
$\End_{\Gamma} ({\E_K}^{r,s})$ is independent of $K$, since the module
${\E_K}^{r,s}$ is tilting. (This follows by a simple induction using
filtrations, and the well-known fact that
$\Hom_{\Gamma}(\Delta(\lambda), \nabla(\mu))$ is $K$ if $\lambda=\mu$
and 0 otherwise.)

Thus, in order to show that \eqref{SWD} holds for any infinite field
$K$, one need only prove that the dimensions of $\tau_K(
\Br_{r,s}^{(n)} )$ and $\End_{\Br_{r,s}^{(n)}}({\E_K}^{r,s})$ are
independent of $K$. In other words, the question amounts to knowing
something about the representation theory of $\Br_{r,s}^{(n)}$.  In
case $n\ge r+s$ the authors have shown, by an argument similar to one
in \cite[Theorem 3.4]{DDH}, that $\tau_K$ is injective in any
characteristic, {\em i.e.}, the algebra $\Br_{r,s}^{(n)}$ acts
faithfully, so in that case the independence of $\dim_K \tau_K(
\Br_{r,s}^{(n)} )$ can be established. Details will appear elsewhere.

Since $S_K(n; r,s) =\rho_K(K\Gamma)$, we have in characteristic zero
for $n \ge r+s$ the equality
\begin{equation} \label{cent}
  S_K(n; r,s) = \End_{\Br_{r,s}^{(n)}}({\E_K}^{r,s}).
\end{equation}
If \eqref{SWD} can be established for all $K$, and for all $n$, then
\eqref{cent} will be true without restriction on $n$ or
characteristic.

\bibliographystyle{amsalpha}

\end{document}